\RequirePackage{fix-cm}
\documentclass[smallextended,envcountsame,final]{svjour3}
\smartqed
\usepackage[arrow, matrix, curve]{xy}
\usepackage{natbib}

\usepackage{array,amsmath,amssymb,amsthm,verbatim,epsfig}
\usepackage{amsfonts}%
\usepackage{graphicx}
\usepackage{enumerate}

\newtheorem{fact}[theorem]{Fact}

\newcommand{\N}{\mathbb{N}}
\newcommand{\Aff}{\mathbb{A}}
\newcommand{\suchThat}{\enspace \vert \enspace}

\newcommand{\Lpos}{\operatorname{L}^{+}}

\newcommand{\Loop}{\operatorname{L}}
\newcommand{\Lp}{\operatorname{L}_{p}}
\newcommand{\Lpp}{\operatorname{L}_{p}^{+}}
\newcommand{\dlim}{\underrightarrow{\lim}}
\newcommand{\plim}{\underleftarrow{\lim}}
\newcommand{\W}{\mathbf{W}}

\newcommand{\Lat}[1]{\mathcal{L}att^{#1}}
\newcommand{\sLat}[1]{\mathcal{L}att^{#1,0}}
\newcommand{\pLat}[1]{\mathcal{L}att_{p}^{#1}}
\newcommand{\spLat}[1]{\mathcal{L}att_{p}^{#1,0}}
\newcommand{\pGrass}{\mathcal{G}rass_{p}}
\newcommand{\Gr}{\mathcal{G}rass}
\newcommand{\Dem}[1]{D_{#1}}
\newcommand{\OCell}[1]{C_{#1}}
\newcommand{\SCell}[1]{\mathcal{C}_{#1}}

\newcommand{\catSp}[1]{(\text{Sp}/S)}
\newcommand{\catIndSch}[1]{(\text{ind-Sch}/#1)}
\newcommand{\catSch}[1]{(\text{Sch}/#1)}
\newcommand{\catAlg}[1]{(#1\text{-Alg})}
\newcommand{\catAlgm}[1]{(#1\text{-Alg(m)})}
\newcommand{\catRingedSp}[1]{(\text{Ringed spaces}/#1)}
\newcommand{\catLocRingedSp}[1]{(\text{Loc. ringed spaces}/#1)}
\newcommand{\catSet}{\text{(Set)}}

 \DeclareMathOperator{\rk}{rk}
 \DeclareMathOperator{\Hom}{Hom}
\DeclareMathOperator{\Spec}{Spec} 
\DeclareMathOperator{\Tor}{Tor} 
\DeclareMathOperator{\id}{id}
\DeclareMathOperator{\X}{X}
\DeclareMathOperator{\diag}{diag}

\DeclareMathOperator{\GlobSect}{\Gamma}

\newcommand{\Sl}{\mathbf{SL}}
\newcommand{\Gl}{\mathbf{GL}}
\newcommand{\Mat}{\mathbf{Mat}}

\newcommand{\cochar}{\check{\X}}
\newcommand{\domcochar}{\check{\X}_{+}}

\newcommand{\Up}{U_{\lambda}}

\newcommand{\Vl}{V_{\lambda}}
\newcommand{\Lattice}[1]{\mathcal{L}_{#1}}

\begin{document}

\setcitestyle{authoryear}

\title{On $p$-adic lattices and Grassmannians}
\journalname{Mathematische Zeitschrift}

\author{Martin Kreidl}
\institute{Martin Kreidl, \email{martin.kreidl@uni-due.de}}

\date{Received: 10/06/2011 / Accepted: 00/00/00}

\maketitle

\begin{abstract}
It is well-known that the coset spaces $G(k((z)))/G(k[[z]])$, for a reductive group $G$ over a field $k$, carry the geometric structure of an inductive limit of projective $k$-schemes.
This $k$-ind-scheme is known as the affine Grassmannian for $G$.
From the point of view of number theory it would be interesting to obtain an analogous geometric interpretation of quotients of the form $\mathcal{G}(\W(k)[1/p])/\mathcal{G}(\W(k))$, where $p$ is a rational prime, $\W$ denotes the ring scheme of $p$-typical Witt vectors, $k$ is a perfect field of characteristic $p$ and $\mathcal{G}$ is a reductive group scheme over $\W(k)$.
The present paper is an attempt to describe which constructions carry over from the function field case to the $p$-adic case, more precisely to the situation of the $p$-adic affine Grassmannian for the special linear group $\mathcal{G}=\Sl_{n}$.
We start with a description of the $R$-valued points of the $p$-adic affine Grassmannian for $\Sl_{n}$ in terms of lattices over $\W(R)$, where $R$ is a perfect $k$-algebra.
In order to obtain a link with geometry we further construct projective $k$-subvarieties of the multigraded Hilbert scheme which map equivariantly to the $p$-adic affine Grassmannian.
The images of these morphisms play the role of Schubert varieties in the $p$-adic setting.
Further, for any reduced $k$-algebra $R$ these morphisms induce bijective maps between the sets of $R$-valued points of the respective open orbits in the multigraded Hilbert scheme and the corresponding Schubert cells of the $p$-adic affine Grassmannian for $\Sl_{n}$.



\end{abstract}

\section{Introduction}

Let $k$ denote an arbitrary field, $n\geq 2$ an integer, and let $\Sl_{n}$ denote the special linear group. The affine Grassmannian for $\Sl_{n}$ over $k$ is the quotient fpqc-sheaf
$$
\Gr: R \mapsto \Sl_{n}(R((z)))/\Sl_{n}(R[[z]])
$$
which maps $k$-algebras to sets. This sheaf is represented by an ind-projective $k$-ind-scheme, as we recall briefly in the following.

\begin{definition}[\citealt{beauville-laszlo, goertz10}]\label{defLatFunc}
Let $R$ be a $k$-algebra. A \emph{lattice} $L\subset R((z))^{n}$ is a finitely generated projective $R[[z]]$-submodule such that $L\otimes_{R[[z]]}R((z)) = R((z))^{n}$. A lattice $L$ is called \emph{special}, if its determinant is trivial, i.e. $\wedge^{n}L = R[[z]]$.
\end{definition}

By $\Lat{n}$ we denote the functor which associates to any $R$-algebra the set of lattices in $R((z))^{n}$, and $\sLat{n}$ denotes the subfunctor of special lattices. Further, let $N$ be any positive integer. Then we denote by $\Lat{n}_{N}$ and $\sLat{n}_{N}$ the respective subfunctors of lattices $L$ with the property that $z^{N}R[[z]]^{n} \subset L \subset z^{-N}R[[z]]^{n}$. The following theorem by Beauville and Laszlo establishes the representability of $\Gr$ as a $k$-ind-scheme.

\begin{theorem}[\citealt{beauville-laszlo}]
For any $k$-algebra $R$, the set $\Gr(R)$ is the ascending union $\cup_{N\in\mathbb{N}}\sLat{n}_{N}(R)$, and the functor $\sLat{n}_{N}$ is represented by a closed subscheme of an ordinary Grassmannian (more precisely, the Grassmannian which parametrizes $nN$-dimensional $k$-linear subspaces in $k^{2nN}$). Hence the functor $\Gr$ is an ascending union of projective $k$-schemes, and therefore it is an ind-projective $k$-ind-scheme.
\end{theorem}

The affine Grassmannian, also for other algebraic groups than $\Sl_{n}$, and its variants such as partial or full flag varieties, are well studied as natural objects within the geometric Langlands program and in the theory of local models for certain Shimura varieties, see e.g. \cite{goertz10}.\\

From the point of view of number theory it is perhaps even more interesting to look at quotients of the form $\Sl_{n}(\mathbb{Q}_{p})/\Sl_{n}(\mathbb{Z}_{p})$, or more generally of the form $\Sl_{n}(\W(k)[1/p])/\Sl_{n}(\W(k))$, where $k$ is a perfect field of positive characteristic $p$ and $\W(k)$ denotes the ring of $p$-typical Witt vectors over $k$.
A structure of an ind-scheme on these quotients would lead to an algebraic model of the Bruhat-Tits building of the group $\Sl_{n}(\W(k)[1/p])$.
Let us refer to this setting as the $p$-adic case in what follows, while by the function field case we mean the situation discussed before.

In \cite{haboush} the author has addressed the problem of endowing the quotient sets $\Sl_{n}(\W(k)[1/p])/\Sl_{n}(\W(k))$ with an algebro-geometric structure, and he has shown that the $p$-adic situation is in this respect significantly more complicated than the function field case.
One source of complication in the $p$-adic case is certainly the simple fact that $\W(R)$, with $R$ any ring, does not carry the structure of an $R$-module. For this reason it is not obvious how to construct an analogue of $\sLat{n}_{N}$ inside an inductive limit of classical Grassmannians over $k$.
Indeed, the $p$-adic case is far less understood, and it is still not clear whether it is possible to put a structure of an ind-scheme, or a related algebraic structure, on the quotient sets $\Sl_{n}(\W(k)[1/p])/\Sl_{n}(\W(k))$.

The present paper is an attempt to investigate how much of our geometric understanding of the affine Grassmannian for $\Sl_{n}$ carries over from the function field case to the $p$-adic case.
To this end we make the following definitions in analogy to the function field case.
Let $k$ be a perfect field of characteristic $p>0$.

\begin{definition}\label{defnPadicGrassmannian}
The $p$-adic affine Grassmannian for $\Sl_{n}$ is the fpqc-sheaf associated to the functor
$$
\pGrass: R \mapsto \Sl_{n}(\W(R)[1/p])/\Sl_{n}(\W(R))
$$
from the category of $k$-algebras to the category of sets.
\end{definition}

In Section \ref{sectLattices} we will obtain a description of the set of $R$-valued points of the $p$-adic affine Grassmannian for $\Sl_{n}$ in terms of lattices in the case when $R$ is a \emph{perfect} $k$-algebra.

\begin{definition}\label{defLatPadic}
Let $R$ be a perfect $k$-algebra. A \emph{lattice} $L\subset \W(R)[1/p]^{n}$ (or simply: a $\W(R)$-lattice of rank $n$) is a finitely generated, projective $\W(R)$-submodule $L\subset \W(R)[1/p]^{n}$ such that $L\otimes_{\W(R)}\W(R)[1/p] = \W(R)[1/p]^{n}$.
Further, a lattice $L\subset \W(R)^{n}$ is called \emph{special} if $\wedge^{n}L = \W(R)$.
\end{definition}

First, we obtain in Subsection \ref{subsectPadicLattices} the following characterization of lattices in $\W(R)[1/p]^{n}$.
A $\W(R)$-submodule $L\subset\W(R)[1/p]^{n}$ is a lattice if and only if it is a free $\W(R)-module$ Zariski-locally on $R$ and if and only if it is a free $\W(R)$-module fpqc-locally on $R$.
This is again analogous to a well known characterization of lattices in the function field case, see \cite{beauville-laszlo} or Subsection \ref{subsectPadicLattices} in the present paper.

The connection between the $p$-adic affine Grassmannian for $\Sl_{n}$ and the notion of lattice of rank $n$ in the $p$-adic setting will then be established by the following theorem, which we are going to prove in Subsection \ref{subsectLatticesGrassmannian}.

\begin{theorem}\label{thmIntro2}
If $R$ is a perfect $k$-algebra, then the set of $R$-valued points of $\pGrass$ is canonically identified with the set of special lattices $L\subset \W(R)[1/p]^{n}$.
\end{theorem}

As already indicated above, our main goal will then be to investigate the $p$-adic affine Grassmannian for $\Sl_{n}$ from the viewpoint of algebraic geometry.
To this end we review the notions of ind-schemes (Subsection \ref{subsectIndSchemes}) as well as Greenberg realizations (Subsection \ref{subsectGreenbergRealizations}) and \emph{localized} Greenberg realizations (Subsection \ref{subsectLocalizedGreenberg}), which were introduced in a similar manner in \cite{haboush}.
Moreover, building on these notions we introduce in Subsection \ref{subsectLoopGroups} the $p$-adic analogues of the (algebraic) loop group, $\Lp\Sl_{n}$, and the positive (algebraic) loop group, $\Lpp\Sl_{n}$.
In terms of $R$-valued points, where $R$ is any $k$-algebra, we have $\Lpp\Sl_{n}(R) = \Sl_{n}(\W(R))$ and $\Lp\Sl_{n}(R) = \Sl_{n}(\W(R)[1/p])$.
Obviously, these functors come with a natural morphism $\Lpp\Sl_{n}\to\Lp\Sl_{n}$.
With these definitions we can state that the $p$-adic affine Grassmannian for $\Sl_{n}$ is the fpqc-sheaf quotient of loop groups $\Lp\Sl_{n}/\Lpp\Sl_{n}$ (Definition \ref{defnRedefPadicGrassmannian}).
Further, to each dominant cocharacter $\lambda$ of the standard maximal torus $T\subset \Sl_{n}$ we associate a $k$-valued point of $\pGrass$ and let $\SCell{\lambda}$ be its orbit for the natural left-action of $\Lpp \Sl_{n}$ on $\pGrass$.
The orbits $\SCell{\lambda}$ play the role of Schubert cells in $\pGrass$.

In order to link these orbits to $k$-schemes and consider their closures in an appropriate setting we recall in Subsection \ref{subsectHilbert} the construction of multigraded Hilbert schemes by Haiman and Sturmfels and introduce in Subsection \ref{subsectLatticeSchemes} the notion of lattice schemes inside an appropriate affine space.
We think of this construction as the analogue of considering lattices as certain \emph{linear} subspaces inside an affine space in the function field case.
In our setting we obtain, in Subsection \ref{subsectionLatticesWittVector}, for each dominant cocharacter $\lambda$, a projective $k$-subvariety $\Dem{\lambda}$ of a multigraded Hilbert scheme, which parametrizes certain lattice schemes and carries a natural $\Lpp\Sl_{n}$-action.
The link to the $p$-adic affine Grassmannian for $\Sl_{n}$ is then established by the following theorem, which is proved in the course of Subsections \ref{subsectMorphism} and \ref{subsectPropertiesOfMorphism}.

\begin{theorem}\label{thmIntro1}
For every dominant cocharacter $\lambda$ of the standard maximal torus $T\subset\Sl_{n}$ there is an $\Lpp\Sl_{n}$-equivariant morphism of fpqc-sheaves $\pi_{\lambda}: \Dem{\lambda}\to \pGrass$ which has the following properties.
Let $\OCell{\lambda}\subset \Dem{\lambda}$ be the open orbit, and let $\SCell{\lambda} \subset \pGrass$ be the Schubert cell corresponding to $\lambda$.
Then $\pi_{\lambda}$ induces bijections $\OCell{\lambda}(R)\simeq\SCell{\lambda}(R)$ for all reduced $k$-algebras $R$.
Moreover, the image under $\pi_{\lambda}$ of $\Dem{\lambda}(k)$ is precisely the union of the sets of $k$-valued points of the Schubert cells $\SCell{\lambda'}$, with $\lambda'\leq\lambda$ for the classical Bruhat-order (whose definition is recalled in Subsection \ref{subsectPropertiesOfMorphism}).
\end{theorem}

The morphisms $\pi_{\lambda}$ are not injective at the level of $k$-valued points.
This means in particular that the varieties $\Dem{\lambda}$ fail to be good analogues of Schubert varieties in the function field case.
Rather, the situation appears similar to constructions in \cite{kreidl-2010}, and thus the varieties $\Dem{\lambda}$ could perhaps be viewed as a sort of Demazure resolution in the $p$-adic setting.
We make this explicit in the simplest non-trivial special case $n=2$ and $\lambda=(1,-1)$.

Finally, in the appendix we collect a couple of general results on fpqc-sheaves and fpqc-sheafification which are used throughout the paper.
Moreover, we discuss briefly the set-theoretical problems which occur when talking about fpqc-sheafifications, and which are often ignored.
Using results from \cite{waterhouse} we check that such complications do not occur in our construction of the $p$-adic affine Grassmannian for $\Sl_{n}$ as an fpqc-sheaf quotient of loop groups.

\section{Lattices over the Witt ring}\label{sectLattices}

Here and for the rest of this paper we denote by $k$ a perfect field of positive characteristic $p$.
The main goal of this section is to describe the $R$-valued points of the $p$-adic affine Grassmannian for $\Sl_{n}$ in terms of lattices if $R$ is a perfect $k$-algebra.

\subsection{Lattices over $\W(R)$ are free $\W(R)$-modules locally on $R$}\label{subsectPadicLattices}

Let us first recall the following result by \cite{beauville-laszlo} from the function field case.

\begin{theorem}\label{thmFunctionLattices}
For an $R[[z]]$-submodule $L\subset R((z))^{n}$ the following four statements are equivalent:
\begin{enumerate}[(1)]
\item The submodule $L$ is a lattice.
\item Zariski-locally on $R$, $L$ is a free $R[[z]]$-submodule of rank $n$ (i.e. there exist $f_{1},\dotsc,f_{r}\in R$ such that $(f_{1},\dotsc,f_{r})=R$ and for all $i$, $L\otimes_{R[[z]]}R_{f_{i}}[[z]]$ is free of rank $n$ and $L\otimes_{R[[z]]}R((z)) = R((z))^{n}$).
\item fpqc-locally on $R$, $L$ is a free $R[[z]]$-submodule of rank $n$ (i.e. there exists a faithfully flat ring homomorphisms $R\to S$ such that $L\otimes_{R[[z]]}S[[z]]$ is free of rank $n$ and $L\otimes_{R[[z]]}R((z)) = R((z))^{n}$).
\item There exists a positive integer $N$ such that $z^{N}R[[z]]^{n} \subset L \subset z^{-N}R[[z]]^{n}$ and $z^{-N}R[[z]]^{n}/L$ is a projective $R$-module.
\end{enumerate}
\end{theorem}

The statement that (1), (2) and (3) in this theorem are equivalent can be rephrased:
The functor $\Lat{n}$, as defined in the introduction, is the Zariski- resp. the fpqc-sheafification of the functor which associates to each $k$-algebra $R$ the set of \emph{free} lattices in $R((z))^{n}$.
Similarly, $\sLat{n}$ is the Zariski- resp. fpqc-sheafification of the functor which associates to each $R$ the set of \emph{free} special lattices in $R((z))^{n}$.
Our goal in this subsection is to obtain a similar result in the Witt vector setting in the case where $R$ is a perfect $k$-algebra.\\

In what follows we build on our definition of the $p$-adic affine Grassmannian for $\Sl_{n}$ and our notion of $p$-adic lattices as in Definitions \ref{defnPadicGrassmannian} and \ref{defLatPadic} in the introduction.
By $\pLat{n}(R)$ we denote the set of lattices of rank $n$ over $\W(R)$, and $\spLat{n}(R)\subset \pLat{n}(R)$ denotes the subset of special lattices.
First we want to see that the assignment $R\mapsto \pLat{n}(R)$ is a functor on the category of perfect $k$-algebras. To this end we prove the following lemma.

\begin{lemma}\label{lemTensor}
Let $R\to S$ be a homomorphism of perfect $k$-algebras, and let $L$ be a flat $\W(R)$-module satisfying $p^{N}\W(R)^{n}\subset L\subset p^{-N}\W(R)^{n}$ for some positive integer $N$. Then we have
\begin{equation*}
\Tor_{1}^{W(R)}(p^{-N}\W(R)^{n}/L, \W(S)) = 0.
\end{equation*}
In particular, this implies that $L\otimes_{\W(R)}\W(S)\subset p^{-N}\W(S)^{n} \subset \W(S)[1/p]^{n}$.
\end{lemma}

\begin{proof}
Set $F_{R}=p^{-N}\W(R)^{n}, F_{S}=p^{-N}\W(S)^{n}=F_{R}\otimes_{\W(R)}\W(S)$ and consider the exact complex
\begin{multline*}
0 \to \Tor_{1}^{W(R)}(F_{R}/L, \W(S)) \to L\otimes_{\W(R)}\W(S) \to \\ \to F_{S} \to F_{S}/L\otimes_{\W(R)}\W(S) \to 0.
\end{multline*}
On the one hand, multiplication by $p^{2N}$ is the zero-map on $F_{R}/L$, and hence, by functoriality, also on $\Tor_{1}^{W(R)}(F_{R}/L, \W(S))$. On the other hand, $(L\xrightarrow{p}L)\otimes \W(S) = L\otimes(\W(S)\xrightarrow{p}\W(S))$ is injective, as $L$ is flat over $\W(R)$. In other words, $p$ acts faithfully and nilpotently on $\Tor_{1}^{W(R)}(F_{R}/L, \W(S))$, which is thus the 0-module.
\end{proof}

Next we observe that for any finitely generated $\W(R)$-submodule $L\subset \W(R)^{n}$ the condition $L\otimes_{\W(R)}\W(R)[1/p] = \W(R)[1/p]^{n}$ is equivalent to the existence of a positive integer $N$ such that $p^{N}\W(R)^{n} \subset L \subset p^{-N}\W(R)^{n}$.
From this and Lemma \ref{lemTensor} we immediately obtain the following fact.

\begin{fact}
The assignment $R\mapsto \pLat{n}(R)$ defines a functor from the category of perfect $k$-algebras to the category of sets.
Namely, to any homomorphism $R\to S$ assign the map
$$
\pLat{n}(R) \to \pLat{n}(S);\quad L\mapsto L\otimes_{\W(R)}\W(S).
$$
The rule $R\mapsto \spLat{n}(R)$ is a subfunctor.
\end{fact}

The rest of this subsection is devoted to the study of the Zariski- resp. fpqc-sheaf properties of $\pLat{n}$ resp. $\spLat{n}$.

\begin{theorem}\label{thmWittLattices}
\begin{enumerate}[(1)]
\item The functor $\pLat{n}$ is the Zariski-sheafification of the functor on the category of perfect $k$-algebras, which associates to any perfect $k$-algebra $R$ the set of \emph{free} lattices of rank $n$ over $\W(R)$.
\item Moreover, $\pLat{n}$ is even an fpqc-sheaf on the category of perfect $k$-algebras. Together with (1) this says that $\pLat{n}$ is also the fpqc-sheafification of the functor which associates to any perfect $k$-algebra $R$ the set of free lattices of rank $n$ over $\W(R)$.
\item The analogous assertions hold if we replace $\pLat{n}$ by $\spLat{n}$ and ``free lattices of rank $n$'' by ``free special lattices of rank $n$''.
\end{enumerate}
\end{theorem}

\begin{proof}
It suffices to prove the first two parts of the theorem, as part (3) follows directly from (1) and (2).
We will check (1) as follows.
Since by definition $L\in\pLat{n}(R)$ is projective and finitely generated as a $\W(R)$-module, it is even finitely presented and (Zariski-)locally free over $\W(R)$.
This means that there exist $p$-typical Witt vectors $f_{1},\dotsc,f_{m}\in \W(R)$ which generate the unit ideal in $\W(R)$ and such that for each $1\leq i \leq m$ the localization $L\otimes_{\W(R)}\W(R)[1/f_{i}]$ is free over $\W(R)[1/f_{i}]$.
Denote by $g_{i}\in R$ the class mod $p$ of $f_{i}$; without loss of generality we can assume that $g_{i} \neq 0$ for all $i$.
Then the $g_{i}$ generate the unit ideal in $R$, and I claim that the $\W(R[1/g_{i}])$-module $L\otimes_{\W(R)}\W(R[1/g_{i}])$ is free for each $i$.
Denote by $[g_{i}]$ the Teichm\"uller representative of $g_{i}$, which is invertible in $\W(R[1/g_{i}])$.
Thus, in $\W(R[1/g_{i}])$ we may consider the product $[g_{i}]^{-1}f_{i} = \alpha$.
Since the class of $\alpha \pmod p$ is $1$ and $R[1/g_{i}]$ is perfect, we see that $\alpha\in 1+p\W(R[1/g_{i}])$ and in particular that $\alpha$ is invertible in $\W(R[1/g_{i}])$.
Hence the same is true for $f_{i}$, and thus $\W(R)[1/f_{i}] \subset \W(R[1/g_{i}])$.
This proves that we may choose $\coprod_{i=1}^{m}\Spec R[1/g_{i}]\to \Spec R$ as a Zariski-covering over which $L$ is a free lattice.

The proof of part (2) requires more work and will be completed before Corollary \ref{corCharLattices} below.
In what follows, $\W_{N}(R) = \W(R)/p^{N}\W(R)$ denotes the ring of $p$-typical Witt vectors of length $N\in\mathbb{N}$ over a perfect ring $R$.

\begin{lemma}\label{lemTruncate}
Let $R\to S$ be a homomorphism of \emph{perfect} rings. Then
$$
\W_{N}(S)\otimes_{\W_{N}(R)}\W_{N}(S) = \W_{N}(S\otimes_{R}S).
$$
\end{lemma}

\begin{proof}
The ring $\W(S\otimes_{R}S)$ carries a natural structure of a $\W(R)$-algebra, and for this structure we have a homomorphism of algebras $\W(S)\otimes_{\W(R)}\W(S)\to \W(S\otimes_{R}S)$.
We will show by induction on $N$ that this map reduces to an isomorphism modulo $p^{N}$ for every $N$, the case $N=1$ being trivial.
Let us set $A := \W(S)\otimes_{\W(R)}\W(S)$ and $B := \W(S\otimes_{R}S)$, and assume that for $N>1$ the induced map $A/(p^{N-1}) \to B/(p^{N-1})$ is an isomorphism.
Then consider the commutative diagram
\begin{equation*}
\hspace{\fill}
\begin{xy}
\xymatrix{
0 \ar[r] & (p^{N-1}A)/(p^{N}) \ar[d]\ar[r] & A/(p^{N}) \ar[d]\ar[r] & A/(p^{N-1}) \ar[r]\ar[d] & 0 \\
0 \ar[r] & (p^{N-1}B)/(p^{N}) \ar[r] & B/(p^{N}) \ar[r] & B/(p^{N-1}) \ar[r] & 0.
}
\end{xy}
\hspace{\fill}
\end{equation*}
As we assume that $S$ is perfect, multiplication by $p^{m}$ maps $\W(S)$ isomorphically onto $p^{m}\W(S)$ for every $m\in\mathbb{N}$, and thus multiplication by $p^{N-1}$ induces isomorphisms
\begin{equation}\label{eqFirstIso}
S\otimes_{R}S \simeq W(S)\otimes_{\W(R)}(p^{N-1}\W(S)/(p^{N})) \simeq (p^{N-1}A)/(p^{N}).
\end{equation}
On the other hand, also $S\otimes_{R}S$ is a perfect ring, whence multiplication by $p$ is an isomorphism from $\W(S\otimes_{R}S)$ to $p\W(S\otimes_{R}S)$, and
\begin{equation}\label{eqSecondIso}
S\otimes_{R}S \simeq B/(p) \simeq (p^{N-1}B)/(p^{N}).
\end{equation}
Using the isomorphisms \eqref{eqFirstIso} and \eqref{eqSecondIso}, the above commutative diagram becomes isomorphic to
\begin{equation*}
\hspace{\fill}
\begin{xy}
\xymatrix{
0 \ar[r] & S\otimes_{R}S \ar^{\id}[d]\ar[r] & A/(p^{N}) \ar[d]\ar[r] & A/(p^{N-1}) \ar[r]\ar[d] & 0 \\
0 \ar[r] & S\otimes_{R}S \ar[r] & B/(p^{N}) \ar[r] & B/(p^{N-1}) \ar[r] & 0.
}
\end{xy}
\hspace{\fill}
\end{equation*}
As we assumed that the right hand vertical arrow is an isomorphism, the 5-lemma implies that the middle vertical map is an isomorphism, too, which proves the induction step.
\end{proof}

\begin{lemma}\label{lemFaithfullyFlat}
Let $R\to S$ be a homomorphism of perfect rings. Then for every positive integer $N$ the following two statements hold:
\begin{enumerate}[(1)]
\item $\W_{N}(R)\to \W_{N}(S)$ is flat if and only if $R\to S$ is flat,
\item $\W_{N}(R)\to \W_{N}(S)$ is faithful if and only if $R\to S$ is faithful.
\end{enumerate}
(A homomorphism of rings is said to be faithful if and only if it induces a surjective morphism at the level spectra.)
\end{lemma}

\begin{proof}
The first part of the lemma is a special case of the local criterion of flatness \citep[Thm.~22.3]{matsumura-1986}, which includes the statement that if $A$ is a ring, $M$ is an $A$-module and $I$ a nilpotent ideal of $A$ for which $I\otimes_{A}M \simeq IM$ holds, then $M$ is flat over $A$ if and only if $M/IM$ is flat over $A/IA$.
We apply this statement to our situation with $A=\W_{N}(R)$, $M=\W_{N}(S)$ and $I=(p)$, and thus we only have to check that $p\W_{N}(S) = (p)\otimes_{\W_{N}(R)}\W_{N}(S)$.
But this follows from the observation that $(p) = p\W_{N}(R) \simeq\W_{N-1}(R)$, whence $(p)\otimes_{\W_{N}(R)}\W_{N}(S) \simeq \W_{N-1}(S) \simeq p\W_{N}(S)$.

To prove the second statment, we just note that for every ring $R$ the reduction mod $p$, $\W_{N}(R)\to R$, induces a bijection between the associated spectra:
$$
\Spec R\xrightarrow{\simeq} \Spec \W_{N}(R).
$$
Namely, since $p$ is nilpotent in $\W_{N}(R)$ it is contained in every prime ideal of $\W_{N}(R)$.
\end{proof}

\begin{lemma}\label{lemProjective}
Let $(A_{i})_{i\in \mathbb{N}}$ be an inverse system of rings, with all the transition homomorphisms $A_{i+1}\to A_{i}$ surjective, and let $\hat{A}$ be its limit.
Let $M$ be a finitely generated $\hat{A}$-module, write $M_{i}:=M\otimes_{\hat{A}}A_{i}$ and assume that $M = \varprojlim M_{i}$.
If all the $M_{i}$ are projective $A_{i}$-modules, then $M$ is a projective $\hat{A}$-module as well.
\end{lemma}

\begin{proof}
Consider a surjective $\hat{A}$-homomorphism $\pi: \hat{A}^{n}\twoheadrightarrow M$.
We shall show that it splits by constructing a system of compatible splittings of the induced maps $\pi_{i}: A_{i}^{n}\twoheadrightarrow M_{i}$.

Of course, the maps $\pi_{i}$ split, since the $M_{i}$ are projective by assumption.
Our strategy will be to construct \emph{compatible} splittings by induction on $i$.
So assume we have a compatible system of splittings $s_{i}: M_{i}\to A_{i}^{n}$ up to a certain index $i$.
Further we set $I_{i+1}=\ker(A_{i+1}\to A_{i})$, $K_{i+1} = \ker(M_{i+1}\to M_{i})$ and $L_{i+1} = \ker(A_{i+1}^{n}\to A_{i}^{n}) = I_{i+1}^{n}$, respectively.
By tensoring the short exact sequence $0\to I_{i+1} \to A_{i+1}\to A_{i} \to 0$ of $A_{i+1}$-modules with $\pi_{i+1}: (A_{i+1})^{n} \to M_{i+1}$ we obtain the following diagram of $A_{i+1}$-modules, with exact rows and all vertical maps surjective:
\begin{equation*}
\hspace{\fill}
\begin{xy}
 \xymatrix{ 0 \ar[r] & L_{i+1} \ar[d]\ar[r] & A_{i+1}^{n} \ar[r]\ar[d] & A_{i}^{n} \ar[r]\ar[d] & 0\\
            0 \ar[r] & K_{i+1} \ar[r] & M_{i+1} \ar[r] & M_{i} \ar[r] & 0. }
\end{xy}
\hspace{\fill}
\end{equation*}
As the tensor product is right exact, the kernel of $\pi_{i+1}$ surjects onto the kernel of $\pi_{i}$. From this and the 5-lemma we get that the vertical linear map $L_{i+1} \to K_{i+1}$ is onto.

By induction, for the $A_{i}$-linear map $\pi_{i}: A_{i}^{n}\to M_{i}$ we already have a splitting $s_{i}$.
By $A_{i+1}$-projectivity of $M_{i+1}$ we may lift the composition $M_{i+1}\to M_{i}\to A_{i}^{n}$ in order to obtain a map $\widetilde{s_{i+1}}: M_{i+1}\to A_{i+1}$, rendering the right square in the following diagram commutative:
\begin{equation*}
\hspace{\fill}
\begin{xy}
 \xymatrix{ 0 \ar[r] & L_{i+1} \ar[d]\ar[r] & A_{i+1}^{n} \ar[r] & A_{i}^{n} \ar[r] & 0 \\
            0 \ar[r] & K_{i+1} \ar[r] & M_{i+1} \ar[r]\ar[u]^{\widetilde{s_{i+1}}} & M_{i} \ar[u]^{s_{i}}\ar[r] & 0. }
\end{xy}
\hspace{\fill}
\end{equation*}
In general, $\widetilde{s_{i+1}}$ will not be a splitting of $\pi_{i+1}$, but it can be properly adjusted: A diagram chase shows that the difference $\delta_{i} := (\pi_{i+1}\circ \widetilde{s_{i+1}} - 1_{M_{i+1}})$ is an $A_{i+1}$-linear map $M_{i+1} \to \ker(M_{i+1}\to M_{i}) = K_{i+1}$.
As $M_{i+1}$ is a projective $A_{i+1}$-module we can lift $\delta_{i}$ to $\Delta_{i}: M_{i+1}\to L_{i+1} \to A_{i+1}^{n}$ (as remarked above $L_{i+1}\to K_{i+1}$ is surjective).
If we set $s_{i+1} := \widetilde{s_{i+1}} - \Delta_{i}$, we get indeed a splitting of $\pi_{i+1}$, producing a commutative square
\begin{equation*}
\hspace{\fill}
\begin{xy}
 \xymatrix{ A_{i+1}^{n} \ar[r] & A_{i}^{n}\\
            M_{i+1} \ar[r]\ar[u]^{s_{i+1}} & M_{i}\ar[u]^{s_{i}}. }
\end{xy}
\hspace{\fill}
\end{equation*}
Inductively applying this construction, we end up with a projective system of splittings, the limit of which is the desired splitting of $\pi$.
\end{proof}

We are now ready to prove that the functor $R\mapsto \pLat{n}(R)$ is a sheaf for the fpqc-topology on the category of perfect rings.
To begin with, note that for any perfect ring $R$ and any finitely generated $\W(R)$-submodule $M\subset \W(R)[1/p]^{n}$ satisfying $p^{N}\W(R)^{n} \subset M \subset p^{-N}\W(R)^{n}$ for some $N$, we have
\begin{equation}\label{eqProjLim}
\varprojlim (M\otimes \W(R)/p^{i}\W(R)) = \varprojlim M/p^{i}M = \varprojlim M/p^{j}\W(R)^{n} = M.
\end{equation}
Here the second equality holds since the respective inverse systems are coinitial, while the third equality follows from the left exactness of the inverse limit over $j$ of the short exact sequence
$$
0 \to M/p^{j}\W(R)^{n} \to p^{-N}\W(R)^{n}/p^{j}\W(R)^{n} \to p^{-N}\W(R)^{n}/M \to 0. 
$$
Since we already know that $\pLat{n}$ is a Zariski-sheaf, it suffices to consider a faithfully flat homomorphism $R\to S$ of perfect rings, and show that the sequence
\begin{equation}
\pLat{n}(R) \to \pLat{n}(S) \rightrightarrows \pLat{n}(S\otimes_{R}S)
\end{equation}
is an equalizer.

(1) $\pLat{n}(R) \to \pLat{n}(S)$ is injective: Take $L,L'\in \pLat{n}(R)$ such that $L\otimes_{\W(R)}\W(S) = L'\otimes_{\W(R)}\W(S)$.
By Lemma \ref{lemFaithfullyFlat} we know that $\W_{N}(R)\to \W_{N}(S)$ is faithfully flat for every $N$, which tells us that $L\otimes_{\W(R)}\W_{N}(R)=L'\otimes_{\W(R)}\W_{N}(R)$.
Using equation \eqref{eqProjLim} this proves $L=L'$.

(2) The subset $\pLat{n}(R) \subset \pLat{n}(S)$ is the equalizer of $\pLat{n}(S) \rightrightarrows \pLat{n}(S\otimes_{R}S)$: Clearly, $\pLat{n}(R)$ is contained in the difference kernel.
To check the converse, we consider the two $\W(S)$-module structures $j_{1}, j_{2}: \W(S)\rightrightarrows\W(S\otimes_{R}S)$, given by $j_{1}: w\mapsto w\otimes 1$ and $j_{2}: w\mapsto 1\otimes w$, respectively, and form the tensor products $L_{1} = L\otimes_{W(S),j_{1}}\W(S\otimes_{R}S)$ and $L_{2} = L\otimes_{W(S),j_{2}}\W(S\otimes_{R}S)$.
Then $L$ is in the difference kernel if and only if $L_{1} = L_{2}$.
For such an $L$ we may conclude, using Lemma \ref{lemTruncate}, that for $i$ large enough one has
\begin{equation}\label{eqDesc1}
\begin{split}
(L\otimes\W_{i}(S))\otimes_{\W_{i}(S),j_{1}}(\W_{i}(S)\otimes_{\W_{i}(R)}\W_{i}(S))  = \\ = (L\otimes\W_{i}(S))\otimes_{\W_{i}(S),j_{2}}(\W_{i}(S)\otimes_{\W_{i}(R)}\W_{i}(S)),
\end{split}
\end{equation}
and similarly
\begin{equation}\label{eqDesc2}
\begin{split}
(L/(p^{i}\W(S))^{n})\otimes_{\W_{i+N}(S),j_{1}}(\W_{i+N}(S)\otimes_{\W_{i+N}(R)}\W_{i+N}(S)) = \\ = (L/(p^{i}\W(S))^{n})\otimes_{\W_{i+N}(S),j_{2}}(\W_{i+N}(S)\otimes_{\W_{i+N}(R)}\W_{i+N}(S)).
\end{split}
\end{equation}
For $i>2N$ we consider the diagram of $\W_{i+N}(S)$-modules
\begin{small}
\begin{equation*}
\begin{xy}
\xymatrix{
 (p^{-N}\W(S)^{n})/(p^{i+N}\W(S)^{n}) \ar[r] & p^{-N}\W_{i}(S)^{n} \ar@{=}[r] & (p^{-N}\W(S)^{n})/(p^{i-N}\W(S)^{n})\\
 L/(p^{i+N}\W(S))^{n} \ar@{->>}[r]\ar@^{(->}[u] & L\otimes_{\W(S)}\W_{i}(S) \ar@{->>}[r]\ar[u] & L/(p^{i-N}\W(S))^{n} \ar@^{(->}[u] & \\
 (p^{N}\W(S)^{n})/(p^{i+N}\W(S)^{n}) \ar@{=}[r]\ar@^{(->}[u] & p^{N}\W_{i}(S)^{n} \ar[r]\ar[u] & (p^{N}\W(S)^{n})/(p^{i-N}\W(S)^{n}). \ar@^{(->}[u]
}
\end{xy}
\end{equation*}
\end{small}
Now \eqref{eqDesc1} and \eqref{eqDesc2} together with Lemma \ref{lemFaithfullyFlat} say that this diagram descends to a diagram of $\W_{i+N}(R)$-modules, i.e.~we obtain a diagram of the form
\begin{small}
\begin{equation*}
\begin{xy}
\xymatrix{
 (p^{-N}\W(R)^{n})/(p^{i+N}\W(R)^{n}) \ar[r] & p^{-N}\W_{i}(R)^{n} \ar@{=}[r] & (p^{-N}\W(R)^{n})/(p^{i-N}\W(R)^{n})\\
 P_{i+N} \ar@{->>}[r]\ar@^{(->}[u] & M_{i} \ar@{->>}[r]\ar[u] & P_{i-N} \ar@^{(->}[u] & \\
 (p^{N}\W(R)^{n})/(p^{i+N}\W(R)^{n}) \ar@{=}[r]\ar@^{(->}[u] & p^{N}\W_{i}(R)^{n} \ar[r]\ar[u] & (p^{N}\W(R)^{n})/(p^{i-N}\W(R)^{n}). \ar@^{(->}[u]
}
\end{xy}
\end{equation*}
\end{small}
We thus have two cofinal systems of $\W(R)$-modules, $(M_{i})$ and $(P_{i})$, whose inverse limit is a $\W(R)$-module $M$. I claim that this is the desired $\W(R)$-lattice.
First observe that for $N$ big enough we have an exact sequence
$$
0 \to p^{N}\W(R)^{n} \hookrightarrow M \to P_{N} \to 0,
$$
as we can see by taking the inverse limit over $i>0$ of the sequence
$$
0 \to p^{N}\W(R)^{n}/p^{i+N}\W(R)^{n} \hookrightarrow P_{i+N} \to P_{i+N}/p^{N}\W(R)^{n} = P_{N} \to 0.
$$
Since $p^{N}\W(R)^{n}$ is finitely generated, and so is $P_{N}$ by faithfully flat descent, $M$ is finitely generated, too.
On the other hand, as the complex
$$
0 \to p^{N}L\otimes_{\W(S)}\W_{i}(S) \to L\otimes_{\W(S)}\W_{i+N}(S) \to L\otimes_{\W(S)}\W_{N}(S) \to 0
$$
is exact, we obtain by faithfully flat descent a short exact sequence
$$
0 \to p^{N}M_{i} \to M_{i+N} \to M_{N} \to 0.
$$
Passing to the inverse limit over $i$ we obtain
$$
0 \to p^{N}M \to M \to M_{N} \to 0,
$$
and thus $M\otimes_{\W(R)}\W_{N}(R) = M_{N}$, which is a projective $\W_{N}(R)$-module by faithfully flat descent: Indeed for a module over any base ring it is equivalent to say that it is (1) finitely generated and projective or (2) of finite presentation and flat.
By definition these properties are satisfied by $L\otimes_{\W(S)}\W_{N}(S)$, and they descend to $M_{N}$ by \cite{ega-iv}, Prop. 2.5.1 and 2.5.2.
Hence we have arrived at a situation where Lemma \ref{lemProjective} applies, proving that $M = \plim (M\otimes_{\W(R)}\W_{N}(R))$ is a $\W(R)$-lattice.
Clearly, $(M\otimes_{\W(R)}\W(S))\otimes_{\W(S)}\W_{N}(S) = M_{N}\otimes_{\W(S)}\W_{N}(S) = L\otimes_{\W(S)}\W_{N}(S)$.
Taking the limit over $N$ we obtain $M\otimes_{\W(R)}\W(S) = L$, which concludes the proof.
\end{proof}

From Theorem \ref{thmWittLattices} we obtain the following corollary in close analogy to Theorem \ref{thmFunctionLattices}.

\begin{corollary}\label{corCharLattices}
Let $R$ be a perfect $k$-algebra and let $L\subset \W(R)[1/p]^{n}$ be a $\W(R)$-submodule. Then the following three statements are equivalent:
\begin{enumerate}[(1)]
\item The submodule $L$ is a lattice.
\item Zariski-locally on $R$, $L$ is a free $\W(R)$-submodule of rank $n$ (i.e. there exist $f_{1},\dotsc,f_{r}\in R$ such that $(f_{1},\dotsc,f_{r})=R$ and for all $i$, $L\otimes_{\W(R)}\W(R_{f_{i}})$ is free of rank $n$ and $L\otimes_{\W(R)}\W(R)[1/p] = \W(R)[1/p]^{n}$).
\item fpqc-locally on $R$, $L$ is a free $\W(R)$-submodule of rank $n$ (i.e. there exists a faithfully flat ring homomorphisms $R\to S$ such that $L\otimes_{\W(R)}\W(S)$ is free of rank $n$ and $L\otimes_{\W(R)}\W(R)[1/p] = \W(R)[1/p]^{n}$).
\end{enumerate}
\end{corollary}

\begin{proof}
This follows immediately from Theorem \ref{thmWittLattices}.
\end{proof}

It is not clear to me whether there is a good translation of condition (4) of Theorem \ref{thmFunctionLattices} to the Witt vector setting. The obvious obstacle is the fact that $\W(R)$ does not carry a structure of an $R$-module.

\subsection{The $p$-adic affine Grassmannian for $\Sl_{n}$ in terms of lattices}\label{subsectLatticesGrassmannian}

For any perfect $k$-algebra $R$, we obtain from Theorem \ref{thmWittLattices} the following characterization of the $R$-valued points of the $p$-adic affine Grassmannian for $\Sl_{n}$ in terms of lattices.

\begin{theorem}\label{thmLatticesGrass}
The fpqc-sheaf $\spLat{n}$ is equal to the restriction of the $p$-adic affine Grassmannian $\pGrass$ to the category of perfect $k$-algebras.
\end{theorem}

\begin{proof}
The presheaf $R\mapsto \Sl_{n}(\W(R)[1/p])/\Sl_{n}(\W(R))$ coincides with the presheaf $R\mapsto \lbrace \text{free special lattices of rank $n$ over $\W(R)$}\rbrace$ on the category of perfect $k$-algebras.
Thus it suffices to prove the following claim: For each presheaf $F$ on the fpqc-site over $k$ the processes of sheafification and restriction to the category of perfect $k$-algebras commute.
Let $R$ be a perfect $k$-algebra and let $\lbrace U_{i}\to \Spec R \rbrace$ be a covering (on the fpqc-site over $k$).
Refining the covering we may assume that the $U_{i}$ are all affine.
For each $i$ denote by $\widetilde{U_{i}}$ the perfection of $U_{i}$.
Then the morphisms $\widetilde{U_{i}} \to \Spec R$ are still jointly surjective and flat, since the affine ring of the perfection is the direct limit of an inductive system
\begin{equation*}
\hspace{\fill}
\begin{xy}
\xymatrix{
R \ar[d]\ar^{\simeq}[r] & R \ar[d]\ar^{\simeq}[r] & \dotsb \ar^{\simeq}[r] & \dlim R \ar[d]\ar@{=}[r] & R \ar[d]\\
U_{i} \ar[r] & U_{i} \ar[r] & \dotsb \ar[r] & \dlim U_{i} \ar@{=}[r] & \widetilde{U_{i}},
}
\end{xy}
\hspace{\fill}
\end{equation*}
where all horizontal maps are the $p$-th power map, and direct limits preserve flatness. Thus we have obtained a refinement of $\lbrace U_{i}\to\Spec R \rbrace$, which is by definition also a covering in the fpqc-site on the category of perfect $k$-algebras, and the claim now follows from Lemma \ref{lemRestrSheafif} in the appendix.
\end{proof}

We conclude this section by describing the well-known Cartan decomposition for the $p$-adic affine Grassmannian $\pGrass$.
Denote by $T$ the standard maximal torus contained in the standard Borel subgroup $B\subset \Sl_{n}$ of upper triangular matrices, and let $\cochar(T)$ and $\domcochar(T)$ be the sets of cocharacters and dominant cocharacters (respectively). We identify $\cochar(T)$ with the subset of $\mathbb{Z}^{n}$ of vectors whose coordinates sum up to 0, and let $\domcochar(T)\subset \cochar(T)$ be the subset of vectors whose coordinates moreover form a decreasing sequence. Further, consider the embedding
$$
\cochar(T) \hookrightarrow \Sl_{n}(W(k)[1/p]); \quad \lambda=(\lambda_{1},\dotsc,\lambda_{n})\mapsto \diag(p^{\lambda_{1}},\dotsc,p^{\lambda_{n}}).
$$

Composing with the natural morphism of functors $\Sl_{n}(W(R)[1/p]) \to \pGrass(R)$ we obtain a map
$$
\Lattice{}: \domcochar(T) \to \pGrass(k);\quad \lambda\mapsto\Lattice{\lambda},
$$
which is injective by the elementary divisors theorem.

\begin{definition}
For each $\lambda \in \domcochar(T)$ the \emph{Schubert cell} $\SCell{\lambda} \subset \pGrass$ is the fpqc-sheafification of the $\Sl_{n}(W(R))$-orbit of $\Lattice{\lambda}\in\pGrass$.
\end{definition}

\begin{theorem}\label{thmCartan}
The $p$-adic affine Grassmannian for $\Sl_{n}$ is, at the level of $k$-valued points, the disjoint union of its Schubert cells $\SCell{\lambda}$, for all $\lambda\in\domcochar(T)$:
$$
\pGrass(k) = \coprod_{\lambda\in\domcochar(T)}\SCell{\lambda}(k).
$$
\end{theorem}

\begin{proof}
This follows from the elementary divisors theorem.
\end{proof}

\begin{remark}
In this section we have defined the $p$-adic affine Grassmannian for $\Sl_{n}$ and its Schubert cells as fpqc-sheafifications of certain functors.
In general the process of fpqc-sheafification involves set-theoretical complications, with the consequence that in certain cases one cannot speak of such sheafifications without making a non-canonical choice of a universe.
In the appendix to this paper we present an argument (Corollary \ref{corGroupQuotientBasicallyBounded}) in order to prove that these issues do not occur in our situation, and thus our notion of a $p$-adic affine Grassmannian for $\Sl_{n}$ is well-defined.
\end{remark}

\section{Greenberg realizations and loop groups}\label{sectGreenberg}

In this section we recall Greenberg's notion of \emph{realization} (in the category of schemes) and introduce a generalization of Greenberg's definition, which we call \emph{localized Greenberg realization}, in the category of ind-schemes.
The idea of considering localized Greenberg realizations is due to \cite{haboush}, though his definition is slightly different from ours.
We then apply our constructions to obtain a definition of \emph{loop groups} associated with linear algebraic loop groups over a discrete valuation ring, where we are particularly interested in the case of the special linear group $\Sl_{n}$ over the ring of $p$-typical Witt vectors $\W(k)$.

\subsection{The language of ind-schemes}\label{subsectIndSchemes}

In what follows we will make extensive use of the language of ind-schemes. Hence, in this subsection we fix our conventions on ind-schemes and briefly discuss their basic properties.

\begin{definition}
Let $S$ be a scheme. An \emph{$S$-space} is a sheaf on the fpqc-site over $S$. An \emph{ind-scheme over $S$} (or simply \emph{$S$-ind-scheme}) is the filtered colimit in the category of $S$-spaces of a system of quasi-compact $S$-schemes. Morphisms of ind-schemes are morphisms of functors. 

\noindent If an $S$-ind-scheme $X$ has the form $X=\dlim_{i\in I} X_{i}$ with all the $X_{i}$ quasi-compact, then we say that $X$ is \emph{represented} by the direct system $(X_{i})_{i\in I}$. By abuse of language we will also simply speak of the $S$-ind-scheme $(X_{i})_{i\in I}$.

\noindent Let $X$ be an $S$-ind-scheme.
By an \emph{$S$-sub-ind-scheme} $Y\subset X$ we mean a subfunctor of $X$ which is itself an $S$-ind-scheme.
An $S$-sub-ind-scheme $Y\subset (X_{i})_{i}$ is called \emph{ind-closed}, if it is represented by a system of closed subschemes $Y_{i}\subset X_{i}$ and the transition morphisms $Y_{i} \to Y_{j}$ are induced by $X_{i} \to X_{j}$.
\end{definition}

Throughout this paper we will assume that the directed index set $I$ is denumerable. In particular, there always exists a (filtered) cofinal subset $I'\subset I$ which can be identified with the natural numbers.
We denote the category of $S$-schemes by $\catSch{S}$, by $\catSp{S}$ we denote the category of $S$-spaces (the morphisms between two $S$-spaces being natural transformations of functors), and by $\catIndSch{S}$ we denote its full subcategory whose objects are $S$-ind-schemes. In other words, we have the following fully faithful functors:
$$
\catSch{S} \hookrightarrow \catIndSch{S} \hookrightarrow \catSp{S}
$$

\begin{remark}
\begin{enumerate}
\item Our definitions of an $S$-space and an $S$-ind-scheme coincide with those given by \cite{beauville-laszlo} in the case where $S$=$\Spec k$ for some field $k$.
\item The existence of colimits in the category of $S$-spaces of direct systems of $S$-schemes is proved in Proposition \ref{propSheafifOfDirectSystem} of the appendix.
\end{enumerate}
\end{remark}

We collect a few easy facts about ind-schemes.

\begin{lemma}\label{lemMorIndSch}
If $T$ is a quasi-compact $S$-scheme and $X$ is an ind-scheme over $S$ which is represented by a direct system of $S$-schemes $(X_{i})$, then $\Hom_{S}(T,X) = \dlim \Hom_{S}(T,X_{i})$.
\end{lemma}

\begin{proof}
Let $\varphi \in \Hom_{S}(T,X)$.
As we prove in the appendix (Proposition \ref{propSheafifOfDirectSystem}), the ind-scheme $X$ is the Zariski-sheafification of the presheaf-direct limit $\dlim X_{i}$.
Thus we can find a Zariski-covering $\lbrace U_{j}; j\in J\rbrace$ of $T$ ($J$ an arbitrary index set) such that $\varphi$ is determined by a family of morphisms of schemes $\varphi_{j}: U_{j} \to X_{i_{j}}$.
Since $T$ is quasi-compact, we may assume that $J$ is finite.
Further, for each pair $(j,j')\in J^{2}$ there is an index $i_{j,j'}$ such that the two morphisms $U_{j}\cap U_{j'} \to X_{i_{j,j'}}$, induced by $\varphi_{j}$ and $\varphi_{j'}$, respectively, coincide.
Now, if we take $m$ to be the maximum of all the $i_{j}$ and the $i_{j,j'}$, the morphisms $U_{j}\to X_{i_{j}}\to X_{m}$ for $j\in J$ glue to a morphism $T\to X_{m}$ which induces $\varphi$.
\end{proof}

Let $X$ and $Y$ be $S$-ind-schemes which are represented by direct systems $(X_{i})$ and $(Y_{i})$ (respectively) of $S$-schemes.
Any morphism of direct systems $(X_{i})\to (Y_{i})$ (i.e. a system of compatible maps $f_{i}: X_{i}\to Y_{j_{i}}$) induces a morphism $f: X\to Y$.
In this case we say that $f$ is represented by the system $(f_{i})$.
Using Lemma \ref{lemMorIndSch} the following converse is easy to deduce.

\begin{lemma}
Let $X$ and $Y$ be $S$-ind-schemes which are represented by direct systems $(X_{i})$ and $(Y_{j})$ (respectively) of $S$-schemes.
Then every morphism $X\to Y$ is represented by a compatible system of maps $f_{i}: X_{i}\to Y_{j_{i}}$.
\end{lemma}

Note that this lemma holds precisely because quasi-compactness of all the $X_{i}$ is built in the definitions.
Moreover, as remarked above, we can always assume that all our index sets are equal to the set of natural numbers, and that compatible systems of maps are of the form $f_{i}: X_{i}\to Y_{i}$ (i.e. preserve the index).

\begin{lemma}[Products]
Let $X,Y,Z$ be $S$-ind-schemes which are represented by direct systems $(X_{i}),(Y_{i}),(Z_{i})$ (respectively) over $S$, and let $X\to Z$ and $Y\to Z$ be morphisms represented by compatible systems of maps $X_{i}\to Z_{i}$ and $Y_{i}\to Z_{i}$.
Then the fiber product (in the category of $S$-spaces) $X\times_{Z}Y$ is an $S$-ind-scheme and is represented by the direct system $(X_{i}\times_{Z_{i}}Y_{i})$.
\end{lemma}

We will make one further technical assumption to simplify our presentation.
\emph{Throughout this paper, all test-schemes which occur will be assumed to be quasi-compact}.
In other words, all functors are considered to be functors on categories of quasi-compact schemes.
This simplification is justified by the fact that an $S$-space is determined by its values on quasi-compact (or even affine) $S$-schemes.
\subsection{Greenberg realizations}\label{subsectGreenbergRealizations}

Our reference for this is \cite{greenberg-1961}, and we stay close to the notation used there.
Let $S$ be a scheme and $\mathbf{R}\to S$ a ring scheme over $S$.
Then $\mathbf{R}$ represents a sheaf of rings on the Zariski-site over $S$, and thus defines a covariant functor
\begin{equation*}
\begin{split}
G_{\mathbf{R}}: \catSch{S} &\to \catRingedSp{\Spec \mathbf{R}(S)}\\
(X,\mathcal{O}_{X}) &\mapsto G_{\mathbf{R}}(X) = (X,\mathcal{O}_{G_{\mathbf{R}}(X)}),
\end{split}
\end{equation*}
where $\mathcal{O}_{G_{\mathbf{R}}(X)}(U) := \mathbf{R}(U)$, the set of $S$-morphisms from $U$ to $\mathbf{R}$.
The ring scheme $\mathbf{R}$ is called a \emph{local} ring scheme, if the functor $G_{\mathbf{R}}$ has values in the category of locally ringed spaces.

\begin{example} Let $\mathbf{R}=\W_{N}$ be the scheme of $p$-typical Witt vectors of length $N$ over $S=\Spec k$, with $0\leq N\leq \infty$.
We claim that $\W_{N}$ is a local ring scheme.
Namely, for any $S$-scheme $X$ the stalk of $G_{\W_{N}}(X)$ at $x\in X$ is given by $\mathcal{O}_{G_{\W_{N}}(X),x} = \dlim \W_{N}(U)$, and $f=(f_{0},f_{1},\dotsc)\in \mathcal{O}_{G_{\W_{N}}(X),x}$ is invertible if and only if $f_{0}\in\mathcal{O}_{X,x}$ is invertible.
The ``only if''-part is trivial, and the ``if''-part can be seen as follows.
Whenever $f_{0}$ is invertible in $\mathcal{O}_{X,x}$, then there exists an open neighbourhood $U$ of $x$ such that $f_{0}$ is invertible in $\mathcal{O}_{X}(U)$.
But then $f$ is invertible in $\W_{N}(U)$ and a fortiori in $\mathcal{O}_{G_{\W_{N}}(X),x}$.
\end{example}

The situation of this example, $\mathbf{R}$ being the scheme of $p$-typical Witt vectors of finite or infinite length over a perfect field $k$, will be the most interesting for us, as we are aiming towards the construction of $p$-adic loop groups.
Another familiar example of a local ring scheme is the scheme of power series in one variable over $k$, i.e.~the scheme $\Aff^{\mathbb{N}}_{k}$ with the property that for each commutative $k$-algebra $A$ we have an identity of rings $A_{k}^{\mathbb{N}}(A) = A[[z]]$ with $z$ a fixed independent variable.

In the following let $\mathbf{R}$ be a local ring scheme over $S$.

\begin{definition}[\citealt{greenberg-1961}]\label{defnGreenbergReal}
Let $X$ be a scheme over the ring $\mathbf{R}(S)$.
A \emph{Greenberg realization of $X$ over $S$} is an $S$-scheme $F_{\mathbf{R}}X$ which represents the functor
$$
Y \mapsto \Hom(G_{\mathbf{R}}(Y),X),
$$
where $\Hom$ is taken in the category of locally ringed spaces over $\mathbf{R}(S)$.
\end{definition}

In the sequel, to simplify notation, we will occasionally drop the index refering to the ring scheme $\mathbf{R}$.

The following proposition and its corollary are purely formal consequences of the universality of representing objects.
However, since they are especially interesting for our applications in the construction of loop groups, we state them explicitly:

\begin{proposition}\label{propFiberProd}
Realizations commute with fiber products.
More precisely, if $X,X',T$ are $\mathbf{R}(S)$-schemes having Greenberg realizations $FX,FX',FT$ over $S$, then $FX\times_{FT}FX'$ is a Greenberg realization over $S$ of $X\times_{T}X'$.
\end{proposition}

\begin{corollary}\label{propGroupRealizations}
Let $X$ be a group scheme over $R(S)$ which has a Greenberg realization $FX$ over $S$. Then $FX$ is a group scheme over $S$.
\end{corollary}

Let us now explicitly describe Greenberg realizations in situations which are of interest to us (as always, we keep in mind the situation where $S=\Spec k$, and the local ring scheme $\mathbf{R}$ is the scheme of $p$-typical Witt vectors of finite or infinite length). Detailed proofs are presented in \cite{greenberg-1961}.

\begin{proposition}[\citealt{greenberg-1961}]\label{propGreenbergRealAffine}
Assume that there is an isomorphism of $S$-schemes
$$
\varphi = (\varphi_{1},\dotsc,\varphi_{N}): \mathbf{R}\xrightarrow{\simeq} \Aff_{S}^{N},
$$
where $0\leq N \leq \infty$, and let us denote by $\GlobSect$ the functor which associates to any scheme its ring of global sections.
Then the Greenberg realization of $\Aff_{\mathbf{R}(S)}^{d}$ is the $S$-scheme $F(\Aff_{\mathbf{R}(S)}^{d}) = (\Aff^{N}_{S})^{d}$ together with the universal arrow $\lambda: GF(\Aff_{\mathbf{R}(S)}^{d}) \to (\Aff_{\mathbf{R}(S)}^{d})$ which is given in terms of global sections by the ring homomorphism
\begin{equation*}
\begin{split}
\lambda^{\#}: \mathbf{R}(S)[T_{1},\dotsc,T_{d}] &\to \GlobSect(S)[t_{1,1},\dotsc,t_{1,N},\dotsc,t_{d,1},\dotsc,t_{d,N}]\\
T_{i} &\mapsto (t_{i,1},\dotsc,t_{i,N}).
\end{split}
\end{equation*}
If $f: \Aff_{\mathbf{R}(S)}^{d} \to \Aff_{\mathbf{R}(S)}^{e}$ is a morphism of $\mathbf{R}(S)$-schemes and $P_{1},\dotsc,P_{e}$ are the polynomials in $\mathbf{R}(S)[T_{1},\dotsc,T_{d}]$ defining $f$, then the morphism $Ff$ between the respective Greenberg realizations is given in terms of global sections by
\begin{equation*}
\begin{split}
\GlobSect(S)[t'_{1,1},\dotsc,t'_{1,N},\dotsc,t'_{e,1},\dotsc,t'_{e,N}] &\to \GlobSect(S)[t_{1,1},\dotsc,t_{1,N},\dotsc,t_{d,1},\dotsc,t_{d,N}]\\
t'_{i,j} &\mapsto \varphi_{j}(\lambda^{\#}(P_{i})).
\end{split}
\end{equation*}
Here, the $t_{i,j}$ are the coordinates on $F(\Aff_{\mathbf{R}(S)}^{d})$, while the $t_{i,j}'$ are the coordinates on $F(\Aff_{\mathbf{R}(S)}^{e})$.
In other words, to calculate the image of $t'_{i,j}$, we have to substitute $T_{l}\mapsto (t_{l,j})_{j}$ in the polynomial $P_{i}$ and then take the $j$-th component of the result under the isomorphism $\varphi$.
\end{proposition}

\begin{proof}
This is proved for finite $N$ in \citep[Prop. 3]{greenberg-1961}. The proof given there also works for $N = \infty$ without modifications.
\end{proof}

\begin{proposition}[\citealt{greenberg-1961}]\label{propAffineRealizations}
Let $\mathbf{R}$ be a local ring scheme over $S$ which is isomorphic to an $N$-dimensional affine space over $S$ (recall that we allow $N=\infty$).
Let moreover $X$ be an affine scheme of finite type over $\mathbf{R}(S)$ having a Greenberg realization by an affine scheme $FX$ over $S$.
Then every closed subscheme of $X$ has a Greenberg realization over $S$ by a closed subscheme of $FX$.
\end{proposition}

\begin{proof}
This is proved in \cite{greenberg-1961}. The crucial point is the observation that we may, by universality of Greenberg realizations, assume that $X$ itself is an affine space over $\mathbf{R}(S)$. In this case we obtain a Greenberg realization of a closed subscheme $Y\subset X$ as follows.
Let $X=\Aff_{\mathbf{R}(S)}^{d}$ and choose a set of defining equations $f_{m}(X_{1},\dotsc,X_{d}) = 0$ for $Y\subset X$.
Each $X_{i}$ can be viewed as a vector of coordinates $X_{i}=(x_{i,0},\dotsc,x_{i,N})$, according to the isomorphism $\mathbf{R}(S)\simeq \Aff_{S}^{N}(S)$.
Plugging these into the equations $f_{m}=0$ yields coordinate-wise equations in the variables $x_{i,j}$, which are the defining equations of a closed subscheme of $FX$. This subscheme is the Greenberg realization $FY\subset FX$ of $Y\subset X$.
\end{proof}

Let us consider for instance the case $\mathbf{R}=\W_{N}$.
Let $X$ be the affine space $\Aff_{\W_{N}(S)}^{d} = \Spec \W_{N}(S)[T_{1},\dotsc,T_{d}]$.
Then a closed subscheme $X\subset \Aff_{\W_{N}(S)}^{d}$ is given by a set of equations (with $I$ an index set)
$$
\lbrace f_{m}(T_{1},\dotsc,T_{d}) = 0 \suchThat m\in I\rbrace.
$$
The equations of the Greenberg realization $FX\subset \Spec S[t_{i,j}]$ are then obtained by plugging the Witt vectors
$$
(t_{i,0},t_{i,1},\dotsc)\in \W_{N}(S[t_{i,0},t_{i,1},\dotsc])
$$
into the equations $f_{m} = 0$.
Thus the components of the Witt vectors
$$
f_{m}(t_{i,0},t_{i,1},\dotsc)\in \W_{N}(S[t_{i,0},t_{i,1}])
$$
for varying $m$ generate the ideal that defines the Greenberg realization $FX$ of $X$.

\subsection{Localized Greenberg realizations}\label{subsectLocalizedGreenberg}

Let $\mathbf{R}$ be a local ring scheme over a \emph{quasi-compact} scheme $S$.
In this subsection we will generalize Greenberg's notion of realization to the situation where $X$ is a scheme over $\mathbf{R}(S)[1/a]$, for $a\in \mathbf{R}(S)$.
Localized Greenberg realizations will be objects in the category of $S$-ind-schemes.
Again, we remind the reader that the situation of interest to us will be the case where $S=\Spec k$ is the spectrum of a perfect field of positive characteristic $p$, $\mathbf{R}=\W$ is the scheme of $p$-typical Witt vectors over $S$, and $a=p$ is a uniformizer.

Observe that the ring $\mathbf{R}(S)[1/a]$ is the colimit of the inductive system of rings
$$
\mathbf{R}(S)\xrightarrow{\cdot a} \mathbf{R}(S)\xrightarrow{\cdot a} \mathbf{R}(S)\xrightarrow{\cdot a} \dotsb.
$$
Assume again that $\mathbf{R}$ is isomorphic as an $S$-scheme to $\Aff_{S}^{N}$.
By Proposition \ref{propGreenbergRealAffine} the affine line over $\mathbf{R}(S)$ can be realized by the affine space $F(\Aff_{\mathbf{R}(S)}^{1}) = \Aff_{S}^{N}$, and by functoriality of Greenberg realization we obtain the inductive system
$$
\Aff_{S}^{N}\xrightarrow{F(\cdot a)} \Aff_{S}^{N}\xrightarrow{F(\cdot a)} \Aff_{S}^{N}\xrightarrow{F(\cdot a)} \dotsb.
$$
If we denote the corresponding $S$-ind-scheme by
$F_{a}\Aff_{\mathbf{R}(S)}^{1}$,
then for any $S$-scheme $Y$ we obtain natural bijections
\begin{equation*}
\begin{split}
&\Hom_{\catIndSch{S}}(Y,F_{a}\Aff_{\mathbf{R}(S)}^{1}) \simeq \dlim (\Aff_{S}^{N}(Y)) =\\ = &\dlim \Hom_{\catLocRingedSp{R(S)}}(G(Y),\Aff_{\mathbf{R}(S)}^{1}) = \dlim \mathbf{R}(Y) = \mathbf{R}(Y)[1/a].
\end{split}
\end{equation*}
In other words, the functor $Y\mapsto \mathbf{R}(Y)[1/a]$ is represented by the $S$-ind-scheme $F_{a}\Aff_{\mathbf{R}(S)}^{1}$.
This motivates the following definition.

\begin{definition}\label{defnLocGreenbergReal}
Let $X$ be an $\mathbf{R}(S)[1/a]$-scheme.
A \emph{localized Greenberg realization} of $X$ over $S$ is an $S$-ind-scheme which represents the functor $Y\mapsto X(\mathbf{R}(Y)[1/a])$ on the category of (quasi-compact) $S$-schemes.
We denote the localized Greenberg realization of $X$ by $F_{a}X$.
\end{definition}

Since the category of ind-schemes has fiber products, and by the universal property of Greenberg realizations, we obtain:

\begin{enumerate}
\item Let $X \to T$ and $X'\to T$ be morphisms of $\mathbf{R}(S)[1/a]$-schemes which admit localized Greenberg realizations $F_{a}X$, $F_{a}X'$ and $F_{a}T$ over $S$.
Then the fiber product $F_{a}X\times_{F_{a}T}F_{a}X'$ is a localized Greenberg realization over $S$ of $X\times_{T}X'$.
\item If a group scheme $X$ over $\mathbf{R}(S)[1/a]$ has a localized Greenberg realization $F_{a}X$ over $S$, then $F_{a}X$ is a group object in the category of ind-schemes over $S$.
\end{enumerate}

Let us gather a few observations which we will use to prove the existence of localized Greenberg realizations in certain cases.
First note that the existence of a localized Greenberg realization of the affine line $\Aff^{1}_{\mathbf{R}(S)}$ is already proved by our remarks before Definition \ref{defnLocGreenbergReal}.
Now let $X$ be any affine scheme of finite type over $\mathbf{R}(S)[1/a]$ and fix a closed immersion $X\subset \Aff_{\mathbf{R}(S)[1/a]}^{d}$.
Let moreover
$$
\varphi: \Aff_{\mathbf{R}(S)[1/a]}^{d}\to \Aff_{\mathbf{R}(S)[1/a]}^{d}
$$
be the automorphism given by $T_{i} \mapsto aT_{i}$ for $i=1,\dotsc,d$.
This yields a diagram of the form
\begin{equation*}
\hspace{\fill}
\begin{xy}
\xymatrix{
\dotsb \ar[r] & \Aff_{\mathbf{R}(S)[1/a]}^{d} \ar^{\varphi}[r] & \Aff_{\mathbf{R}(S)[1/a]}^{d} \ar^{\varphi}[r] & \cdots \\
\dotsb \ar[r] & \varphi^{n}(X) \ar[u]\ar[r] & \varphi^{n+1}(X) \ar[u]\ar[r] & \cdots,
}
\end{xy}
\hspace{\fill}
\end{equation*}
where all the horizontal maps are isomorphisms of $\mathbf{R}(S)[1/a]$-schemes.
We define $X_{n}$ to be the schematic closure of $\varphi^{n}(X)\hookrightarrow \Aff_{\mathbf{R}(S)[1/a]}^{d} \hookrightarrow \Aff_{\mathbf{R}(S)}^{d}$, which determines an $\mathbf{R}(S)$-ind-scheme $(X_{n})_{n}$.
In the sequel we write for any $\mathbf{R}(S)$-scheme $Y$:
$$
Y[1/a] := Y\times_{\Spec \mathbf{R}(S)}\Spec \mathbf{R}(S)[1/a].
$$
With this notation we have $X_{n}[1/a]\simeq \varphi^{n}(X)\simeq_{\varphi^{-n}} X$ for all $n\in \mathbb{N}$.

\begin{lemma}\label{lemIndSchemeRepr}
The $\mathbf{R}(S)$-ind-scheme $(X_{n})_{n}$ represents the functor
$$L: Y\mapsto \Hom_{\mathbf{R}(S)[1/a]}(Y[1/a],X)$$
on the category of (quasi-compact) $\mathbf{R}(S)$-schemes.
\end{lemma}

\begin{proof}
A morphism of functors $\psi_{n}: X_{n} \to L$ is given by the functorial map 
\begin{multline*}
X_{n}(Y) = \Hom_{\mathbf{R}(S)}(Y,X_{n}) \to \Hom_{\mathbf{R}(S)[1/a]}(Y[1/a],X_{n}[1/a]) \\ 
\simeq_{\varphi^{-n}} \Hom_{\mathbf{R}(S)[1/a]}(Y[1/a],X).
\end{multline*}
Obviously the morphisms $\psi_{n}, n\in\mathbb{N},$ are compatible, so we obtain a morphism of functors $\psi: (X_{n})_{n}\to L$.
Since every $Y[1/a]$-valued point $P$ of $X$ is given by a $d$-tuple $\mathbf{p}$ in
$
\GlobSect(Y[1/a])^{d} = (\GlobSect(Y)\otimes_{\mathbf{R}(S)}\mathbf{R}(S)[1/a])^{d},
$
where $\GlobSect$ denotes the functor of global sections, there exists some $n\in\mathbb{N}$ such that $a^{n}\cdot \mathbf{p}\in \GlobSect(Y)^{d}$. Thus $\varphi^{n}(P)$ extends to a $Y$-valued point of $X_{n}$, which shows that $\psi(Y)$ is surjective for each scheme $Y$ over $\mathbf{R}(S)$.
To check injectivity, take $P,Q\in X_{n}(Y)$ such that $P$ and $Q$ have the same image in $L(Y)$.
This means in particular that the corresponding morphisms
$P',Q': Y[1/a]\to X_{n}[1/a]=\varphi^{n}(X)$ are equal, and consequently the respective $\mathbf{R}(S)$-morphisms $P'',Q'':Y[1/a] \to Y \to X_{n}$ are equal.
But both $P$ and $Q$ are given by $d$-tuples $\mathbf{p},\mathbf{q}$ of sections in $\GlobSect(Y)$, and for these the equality $P''=Q''$ says that there exists an $m\in \mathbb{N}$ such that $a^{m}\mathbf{p}=a^{m}\mathbf{q}$.
This means that the compositions
$$
Y\xrightarrow{P,Q}X_{n}\xrightarrow{\varphi^{m}}X_{n+m}
$$
coincide, whence a fortiori $P$ and $Q$ coincide as elements of $(X_{n})_{n}(Y)$.
\end{proof}

It is now straight forward to construct localized Greenberg realizations for affine $\mathbf{R}(S)[1/a]$-schemes of finite type.

\begin{proposition}
Let $X$ be an affine scheme of finite type over $\mathbf{R}(S)[1/a]$, and assume that $\mathbf{R}$ is isomorphic as an $S$-scheme to some affine space over $S$.
Then there exists an $S$-ind-scheme which represents the functor $Y\mapsto X(\mathbf{R}(Y)[1/a])$ on the category of (quasi-compact) $S$-schemes.
\end{proposition}

\begin{proof}
Fix a closed immersion $X\subset \Aff_{\mathbf{R}(S)[1/a]}^{d}$ and let $(X_{n})_{n}$ be as above.
Now apply Greenberg realization to the $\mathbf{R}(S)$-schemes $X_{n}$ and their transition maps.
I claim that the resulting $S$-ind-scheme $(FX_{n})_{n}$ has the desired form.
Indeed, we have
\begin{multline*}
\Hom(Y,(FX_{n})_{n}) = \dlim\Hom(Y,FX_{n}) =\\= \dlim\Hom_{\mathbf{R}(S)}(\mathbf{R}(Y),X_{n}) = \Hom(\mathbf{R}(Y)[1/a],X),
\end{multline*}
where the second equality is by the definition of Greenberg realization, and the third one follows from Lemma \ref{lemIndSchemeRepr}.
\end{proof}

\begin{example}
Let us illustrate this in our standard situation of $p$-typical Witt vectors of infinite length over a (perfect) field $k$.
Let $X=\Aff_{\W(k)[1/p]}^{d}$.
Then the $k$-ind-scheme which is the localized Greenberg realization of $X$ is given (up to isomorphism) by the inductive system
$$
\Spec k[x_{i,j}; i=1,\dotsc,d; j\in\mathbb{N}] \xrightarrow{\cdot p} \Spec k[x_{i,j}; i=1,\dotsc,d; j\in\mathbb{N}] \xrightarrow{\cdot p} \dotsc,
$$
where the transition maps $\cdot p$ are defined by $x_{i,j}\mapsto x_{i,j-1}^{p}$ (for $j=1,\dotsc,\infty$) and $x_{i,0}\mapsto 0$.
\end{example}

\subsection{Construction of generalized and $p$-adic loop groups}\label{subsectLoopGroups}

In this subsection we consider the following situation.
Let $\mathbf{D}$ be a local ring scheme over a field $k$ such that $D=\mathbf{D}(k)$ is a discrete valuation ring with uniformizer $u\in D$.
Moreover we assume that $\mathbf{D}$ is isomorphic to $\Aff^{\mathbb{N}}_{k}$ as a scheme over $k$.
Typical special cases are:
\begin{enumerate}
\item The ring scheme of power series in one variable over $k$, and
\item the ring scheme of $p$-typical Witt vectors over a perfect field $k$ of positive characteristic $p$.
\end{enumerate}
By $K$ we denote the field of fractions of $D$.

\begin{definition}
Let $X$ be a scheme over $\Spec K$.
The functor from the category of $k$-algebras to the category of sets,
$$
\Loop X: R \mapsto X(\mathbf{D}(R)[1/u])
$$
will be called the (generalized) loop space associated to $X$. Moreover, if $X$ is a scheme over $\Spec D$, then we call
$$
\Lpos X: R \mapsto X(\mathbf{D}(R))
$$
the (generalized) positive loop space of $X$. By abuse of notation we also write $\Loop X = \Loop (X_{K})$ for a $D$-scheme $X$.
\end{definition}

Obviously, there is a canonical morphism of functors $\Lpos X \to \Loop X$.
If in addition $X=G$ is a group scheme over $D$, then we call $\Loop G$ and $\Lpos G$ the \emph{(generalized) loop group} and the \emph{(generalized) positive loop group}, respectively, associated to $G$.\\

Note that if $\mathbf{D}$ is the $k$-scheme of power series in one variable over $k$, we recover the usual notions of (formal) loop space, loop group etc., as described in \cite{beauville-laszlo}, \cite{pappas-rapoport-2008} and others.
The following proposition is an immediate consequence of our discussion on Greenberg realizations.

\begin{proposition}
If $X$ is an affine scheme of finite type over $D$, then the functor $\Lpos X$ is representable by an affine scheme over $k$, namely the Greenberg realization over $k$ of $X$.
If $X$ is affine and of finite type over $K$, then $\Loop X$ is representable by the localized Greenberg realization over $k$ of $X$.
\end{proposition}

In fact, in all situations that we are going to consider, the affine scheme $X$ comes together with an embedding into some affine space, $X\subset \Aff^{d}_{D}$.
With respect to this embedding, the construction of the localized Greenberg realization $\Loop X$, as described in the preceeding subsection, produces an explicit direct system $(FX_{i})_{i\in \mathbb{N}}$ of $k$-schemes which represents $\Loop X$.
Explicitly, the scheme in the $i$-th step of this direct system parametrizes the $K$-points of $X$ whose coordinates, with respect to the embedding $X\subset \Aff^{d}_{D}$, have poles of order at most $i$ (i.e., belong to $u^{-i}D$).

Specializing the constructions of this subsection to the case where $k$ is a perfect field of positive characteristic $p$, $D=\W(k)$ and $K=\W(k)[1/p]$, we obtain the following objects in analogy to the function field case \citep[cf.][]{beauville-laszlo}.

\begin{definition}
The \emph{$p$-adic loop group} associated with $\Sl_{n}$ over $\W(k)$ is the $k$-ind scheme $\Loop\Sl_{n}$ representing the functor $R\mapsto\Sl_{n}(\W(R)[1/p])$ on the category of $k$-algebras.
Further, the \emph{positive $p$-adic loop group} is the $k$-scheme $\Lpos\Sl_{n}$ which represents the functor $R\mapsto\Sl_{n}(W(R))$.
\end{definition}

To indicate that we are working in the $p$-adic setting, in what follows we will write $\Lp\Sl_{n}$ (resp. $\Lpp\Sl_{n}$) for the (positive) $p$-adic loop groop.
With this notation we may rephrase the definition of $\pGrass$, Definition \ref{defnPadicGrassmannian}.

\begin{definition}\label{defnRedefPadicGrassmannian}
The $p$-adic affine Grassmannian for $\Sl_{n}$ is the fpqc-sheaf quotient $\Lp\Sl_{n}/\Lpp\Sl_{n}$.
\end{definition}

In the next section we will also encounter the $p$-adic loop group associated with $\Gl_{n}$, $\Lp\Gl_{n}$, as well as the $p$-adic loop space $\Lp\Mat_{n}$ associated with the $\W(k)$-scheme $\Mat_{n}$ of $n\times n$-matrices over $\W(k)$.

\section{Hilbert schemes and lattice schemes}\label{sectHilbert}

\subsection{The multigraded Hilbert scheme of Haiman and Sturmfels}\label{subsectHilbert}

We first recall a result by \cite{haiman-2002} on the representability of the multigraded Hilbert functor.

Let $R$ be any ring, and let $\Aff_{R}^{n} = \Spec R[x_{1},\dotsc,x_{n}]$ be the $n$-dimensional affine space over $R$, and identify $u\in \N^{n}$ with the monomial $x_{1}^{u_{1}}\dotsb x_{n}^{u_{n}}$. Then a multigrading of $R[x_{1},\dotsc,x_{n}]$ by a semigroup $A$ is given by a semigroup homomorphism $\deg: \N^{n}\to A$. This induces a decomposition
\begin{equation*}
 R[x_{1},\dotsc,x_{n}] = \oplus_{a\in A} R[x_{1},\dotsc,x_{n}]_{a},
\end{equation*}
where $R[x_{1},\dotsc,x_{n}]_{a}$ is the $R$-span of the monomials of degree $a$.

A homogeneous ideal $I \subset R[x_{1},\dotsc,x_{n}]$ is called an \emph{admissible ideal}, if for each $a\in A$ the graded piece $(R[x_{1},\dotsc,x_{n}]/I)_{a}$ is a locally free module of constant finite rank on $\Spec R$. Thus every admissible ideal $I \subset R[x_{1},\dotsc,x_{n}]$ has a well-defined \emph{Hilbert function}, given by
\begin{equation*}
 h_{I}: A \to \N, \quad a\mapsto \rk(R[x_{1},\dotsc,x_{n}]/I)_{a}.
\end{equation*}
A closed $R$-subscheme $V\subset \Spec R[x_{1},\dotsc,x_{n}]$ which is defined by an admissible ideal will also be called admissible, and by the Hilbert function of $V$ we mean the Hilbert function of its defining ideal.

Let $h: A\to \N$ be any function vanishing outside $\deg(\N^{n})\subset A$, and define the \emph{Hilbert functor} $\mathcal{H}_{R}^{h}$ from the category of $R$-algebras to sets by
\begin{multline*}
 \mathcal{H}_{R}^{h}(S) = \lbrace \text{admissible ideals } I \subset S[x_{1},\dotsc,x_{n}] \suchThat \\ \rk(S[x_{1},\dotsc,x_{n}]/I)_{a} = h(a) \text{ for all } a\in A \rbrace.
\end{multline*}

\begin{theorem}[Haiman, Sturmfels]\label{thmGradedHilbertscheme}
There exists a quasiprojective scheme $H_{R}^{h}$ over $R$ which represents the functor $\mathcal{H}_{R}^{h}$. If the grading of $R[x_{1},\dotsc,x_{n}]$ is positive, i. e. $1$ is the only monomial with degree $0$, then this scheme is even projective over $R$.
\end{theorem}

The scheme $H_{R}^{h}$ is called the \emph{multigraded Hilbert scheme} for the Hilbert function $h$. If no Hilbert function $h$ is specified, we will use the term ``multigraded Hilbert scheme'' to refer to the disjoint union of the $H_{R}^{h}$ for all Hilbert functions $h$. We denote this disjoint union by $H_{R}$, or simply by $H$ if the ring $R$ is clear from the context.

\subsection{Lattice schemes in general}\label{subsectLatticeSchemes}

For any ring scheme $\mathbf{X}$ over $R$ we have the obvious notion of an $\mathbf{X}$-module scheme over $R$. In particular, we have the free $\mathbf{X}$-module scheme of rank $n$, denoted $\mathbf{X}^{n}$. An $\mathbf{X}$-submodule scheme of an $\mathbf{X}$-module scheme $M$ is a closed $R$-subscheme of $M$ which is stable under the morphisms defining the module operations on $M$.
This means that a closed $\mathbf{X}$-subscheme $V\subset M$ is an $\mathbf{X}$-submodule scheme if the following diagrams exist,
\begin{equation*}
\hspace{\fill}
\begin{xy}
\xymatrix{ M \times M \ar^{add.}[rr] & & M \\
 V\times V \ar@{.>}[rr]\ar@^{(->}[u] & & V, \ar@^{(->}[u]
}
\end{xy}\qquad
\begin{xy}
\xymatrix{ \mathbf{X} \times M \ar^{mult.}[rr] & & M \\
 \mathbf{X}\times V \ar@{.>}[rr]\ar@^{(->}[u] & & V, \ar@^{(->}[u]
}
\end{xy}
\hspace{\fill}
\end{equation*}
and analogous diagrams exist for the zero-section and additive inverses.

In what follows, we always assume that $\mathbf{X}$ is a ring scheme which is isomorphic as an $R$-scheme to $\Aff_{R}^{N}$ ($0\leq N < \infty$). Let us furthermore fix a grading over $R$ of the structure sheaf of $\mathbf{X}\simeq\Aff_{R}^{N}$ so that the ring operations on $\mathbf{X}$ are defined by graded homomorphisms on the structure sheaf. Then also the structure sheaf of $\mathbf{X}^{n}$ is graded. We call a submodule scheme in $\mathbf{X}^{n}$ a \emph{lattice-scheme} if its defining ideal is admissible.

\begin{proposition}
The set of lattice schemes in $\mathbf{X}^{n}$ with a given Hilbert function $h$ is parametrized by a closed subscheme $Z$ of the multigraded Hilbert scheme of $\mathbf{X}^{n}$ over $R$. The $R$-scheme $Z$ is quasi-projective, and it is projective over $R$ if the grading of $\mathbf{X}$ is positive.
\end{proposition}

\begin{proof}
Let $H\to \Spec R$ be the multigraded Hilbert scheme of $\mathbf{X}^{n}$ and let $U\to H$ be the universal family. We have to show that there exists a closed subscheme $Z \subset H$ such that for any morphism $Y\to H$, $V = Y\times_{H}U\subset Y\times_{\Spec R}\mathbf{X}^{n}$ is a submodule scheme if and only if $Y\to H$ factors through $Z \subset H$. It suffices to check this locally on $H$, i.e. for an affine open subscheme $H'=\Spec S \subset H$. Then also $U' := H'\times_{H}U$ is affine, and $U'$ is defined by an ideal $I\subset S[x_{i,j}\suchThat i=1,\dotsc,n; j=0,1,\dotsc,N-1]$ such that the quotient $S[x_{i,j}]/I$ is $S$-locally free. Now for any morphism $Y'= \Spec S' \to H'$ the condition that $V' = Y'\times_{H'}U' \subset U'$ be stable under addition on $\mathbf{X}^{n}$ translates into the condition that its defining ideal $I$ maps to $0$ under the cohomorphism of addition. Analogous vanishing conditions hold for scalar multiplication, units and additive inverses. Since $S[x_{i,j}
]/I$ is locally free over $S$, these vanishing conditions can be expressed by equations with coefficients in $S$, which then define a closed subscheme $Z' \subset H' = \Spec S$. By construction, $V'$ is stable under the module operations if and only if $Y'\to H'$ factors throuth $Z'$. By gluing all the $Z'\subset H$ we obtain the closed subscheme $Z\subset H$ with the desired universal property.
\end{proof}

\begin{proposition}[Group actions on $H$]\label{propGroupActions}
Let $G/\Spec R$ be an algebraic group acting algebraically on $\mathbf{X}^{n}$, and assume that this action respects the grading on the structure sheaf of $\mathbf{X}^{n}$. Then $G$ acts equivariantly on the Hilbert scheme $H$ of $\mathbf{X}^{n}$ and its associated universal family. If furthermore the action of $G$ on $\mathbf{X}^{n}$ is by automorphisms of $\mathbf{X}$-module schemes, then the action of $G$ on the universal family over $H$ restricts to an equivariant action on the universal family over $Z$.
\end{proposition}

\begin{proof}
This is a formal consequence of the universal properties of $H$ and $Z$ and the fact that the action of $G$ on $\mathbf{X}^{n}$ is algebraic, i.e. functorial.
\end{proof}

\subsection{Lattice schemes in the Witt vector setting}\label{subsectionLatticesWittVector}

In what follows, $k$ denotes a perfect field of positive characteristic $p$, and $R$ denotes a $k$-algebra.
Let us specialize the above discussion to the case where
$$
\mathbf{X}=\W_{2N}=\Spec k[\alpha_{-N},\dotsc,\alpha_{N-1}]
$$
is the scheme of Witt-vectors over $S=\Spec k$ of length $2N<\infty$.

Consider further the Greenberg realization over $k$ of the $n$-dimensional affine space $\Aff_{\W_{2N}(k)}^{n}$, which is an affine space over $k$ and carries, by functoriality of Greenberg realization, an obvious structure of $\W_{2N}$-module scheme:
$$
\W_{2N}^{n} \simeq \Aff_{k}^{2N\times n} = \Spec k[x_{i,j} \suchThat i=1,\dotsc,n; j=-N,\dotsc,N-1].
$$
From Witt vector arithmetics \citep[chap.~II, \S 6]{serre} it follows that the morphisms defining the module operations on $\W_{2N}^{n}$ are defined by \emph{graded} homomorphisms of the respective affine rings if we set $\deg(1) = 0, \deg \alpha_{j} = p^{j}\text{ and } \deg x_{i,j} = p^{j}$. This way $\deg$, for each of the respective coordinate rings, is a semigroup homomorphism with values in $p^{-N}\mathbb{N}$. Note that the standard grading with $\deg x_{i,j} = 1$ for all $i,j$ is not respected by the module-operations on $\W_{2N}^{n}$ and is thus not suited for our constructions.

Let us now discuss lattice schemes inside $\W_{2N}^{n}$.
For each dominant cocharacter $\lambda\in\domcochar(T)$ and $N$ such that $N\geq\lambda_{1}\geq\dotsc\geq\lambda_{n}\geq -N$, consider the subscheme $\Vl^{(N)} \subset \W_{2N}^{n}$ defined by the ideal
\begin{equation}\label{eqDefnIl}
I_{\lambda}^{(N)} = ( x_{1,-N},\dotsc,x_{1,\lambda_{1}-1},\dotsc,x_{n,-N},\dotsc,x_{n,\lambda_{n}-1} ) .
\end{equation}
The subscheme $\Vl^{(N)}$ is a lattice scheme, and we denote by $\OCell{\lambda}^{(N)}$ its orbit in $H$ under the action of the linear $k$-group $\Lpp \Sl_{n}$, and by $\Dem{\lambda}^{(N)}$ its orbit-closure.
Theorem \ref{thmGradedHilbertscheme} asserts in particular that $\Dem{\lambda}^{(N)}$ is a projective $k$-variety, which contains $\OCell{\lambda}^{(N)}$ as an open subvariety. 

Clearly, $\Vl^{(N)}$ as well as $\Dem{\lambda}^{(N)}\subset H$ depend on our particular choice of $N$.
However, this choice does not really matter: If we choose $N'>N$, then $\Vl^{(N')}$ and $\Dem{\lambda}^{(N')}$ lie in a Hilbert scheme for a different affine space, but the functorial map
\begin{equation}\label{eqDefIl}
\Dem{\lambda}^{(N)} \to \Dem{\lambda}^{(N')};\quad V \mapsto V\times\Aff_{k}^{N'-N}
\end{equation}
defines a natural $\Lpp\Sl_{n}$-equivariant isomorphism $\Dem{\lambda}^{(N)} \simeq \Dem{\lambda}^{(N')}$ which takes $\Vl^{(N)}$ to $\Vl^{(N')}$.
In what follows, we will thus drop the upper index $(N)$ and write $\OCell{\lambda} \subset \Dem{\lambda}$.
Further, let $\Up^{(N)}\subset\W_{2N}^{n}\times_{\Spec k}\Dem{\lambda}$ be the universal family over $\Dem{\lambda}$ and let $\mathbf{\Up^{(N)}}$ be the preimage of $\Up^{(N)}$ under the natural projection
\begin{equation*}
\W^{n} \times_{\Spec k}\Dem{\lambda} \to \W_{2N}^{n}\times_{\Spec k}\Dem{\lambda}.
\end{equation*}
Then the morphism $\mathbf{\Up^{(N)}} \to \Dem{\lambda}$ is by construction equivariant for the action of $\Lpp\Sl_{n}$. Let us also note, that as an abstract $k$-scheme $\mathbf{\Up^{(N)}}$ is independent of $N$, but its embedding into $\W^{n}\times_{\Spec k}\Dem{\lambda}$ is not. From these constructions we derive the following proposition.

\begin{proposition}\label{propMorKPoints}
For each dominant cocharacter $\lambda\in\domcochar(T)$ there exists a projective $k$-variety $\Dem{\lambda}$ together with a $\Dem{\lambda}$-sub-ind-scheme
$$
\mathcal{U}_{\lambda} \subset \Lp\W(k)^{n}\times_{\Spec k}\Dem{\lambda}
$$
which is invariant for the action of $\Lpp\Sl_{n}$, and such that the map
\begin{equation}
\Dem{\lambda}(k) \to \pGrass(k);\quad V \mapsto (\mathcal{U}_{\lambda}\times_{\Dem{\lambda}}V)(k)
\end{equation}
is well-defined and $\Lp\Sl_{n}(k)$ equivariant.
Under this map, the image of $\Vl\in\Dem{\lambda}$ is equal to the lattice $\Lattice{\lambda}\subset\W(k)[1/p]^{n}$.
\end{proposition}

\begin{proof}
For each positive integer $N$ which satisfies
\begin{equation}\label{eqIneqN}
N\geq\lambda_{1}\geq\dotsc\geq\lambda_{n}\geq -N
\end{equation}
we define $\mathbf{\Up^{(N)}}$ as in the paragraph before the statement of the proposition.
For those finitely many $N$, for which \eqref{eqIneqN} does not hold, we set $\mathbf{\Up^{(N)}}$ equal to the empty scheme $\emptyset$.
Further, we consider an inductive system representing the localized Greenberg realization over $k$ of $\Aff_{\W(k)[1/p]}^{n} = \W(k)[1/p]^{n}$:
\begin{equation}\label{eqIndSystW}
\W^{n} \xrightarrow{\cdot p} \W^{n} \xrightarrow{\cdot p} \cdots,
\end{equation}
where $\cdot p$ denotes the morphism of $k$-schemes which arises from multiplication by $p$ via Greenberg realization.
Now base change from $\Spec k$ to $\Dem{\lambda}$ applied to \eqref{eqIndSystW} yields an inductive systems of $\Dem{\lambda}$-schemes
\begin{equation*}
\begin{xy}
\xymatrix{
\mathbf{\Up^{(1)}} \ar@^{(->}[d] &
\cdots &
\mathbf{\Up^{(N)}} \ar@^{(->}[d] &
\mathbf{\Up^{(N+1)}} \ar@^{(->}[d] & \cdots\\
W^{n}\times\Dem{\lambda} \ar[r]^{\cdot p} &
\cdots \ar[r]^{\cdot p} &
W^{n}\times\Dem{\lambda} \ar[r]^{\cdot p} &
W^{n}\times\Dem{\lambda} \ar[r]^{\cdot p} & \cdots,
}
\end{xy}
\end{equation*}
and each of the vertical morphisms in this system is equivariant for the action of $\Lp\Sl_{n}$ by construction.
In fact, the horizontal morphisms restrict to an inductive system
\begin{equation}\label{eqIndSystU}
\mathcal{U}_{\lambda}: \quad \mathbf{\Up^{(1)}} \xrightarrow{\cdot p} \cdots \xrightarrow{\cdot p} \mathbf{\Up^{(N)}} \xrightarrow{\cdot p} \mathbf{\Up^{(N+1)}} \xrightarrow{\cdot p} \cdots,
\end{equation}
which can be seen as follows.
Obviously, over the point $\Vl$ the horizontal morphisms (multiplication by $p$) on $\W^{n}$ restrict to morphisms between the fibers in the respective $\mathbf{U_{\lambda}^{(N)}}$, by definition of the ideals $I_{\lambda}^{(N)}$.
By $\Lp\Sl_{n}$ equivariance, we can extend this observation to the open orbit $\OCell{\lambda}$, and since $\OCell{\lambda}\subset\Dem{\lambda}$ is dense and $\mathbf{U}^{(N)}_{\lambda}\subset\W^{n}\times\Dem{\lambda}$ is closed for all $N$, this holds on all of $\Dem{\lambda}$.
Thus we obtain an $\Lp\Sl_{n}$-equivariant morphism of $k$-ind-schemes
$$
\mathcal{U}_{\lambda} \hookrightarrow \Lp\W(k)[1/p]^{n}\times_{\Spec k}\Dem{\lambda}.
$$
Next we consider the fiber $\mathcal{F}_{\lambda}$ of the sub-ind-scheme $\mathcal{U}_{\lambda} \subset \Lp\W(k)[1/p]^{n}\times\Dem{\lambda}$ over $\Vl$.
Assume that $N$ satisfies the inequalities \eqref{eqIneqN}.
Then, since $k$ is perfect, each $k$-valued point of $\mathcal{F}_{\lambda}$ is given by a morphism
\begin{equation*}
\Spec k \to F_{\lambda} := \mathbf{\Up^{(N)}}\times_{\Dem{\lambda}}\lbrace\Vl\rbrace = \Spec (\GlobSect(\W^{n})/(I_{\lambda}^{(N)})) \subset \W^{n},
\end{equation*}
where $\GlobSect(\W^{n}) = \Spec k[x_{i,j}\suchThat i=1,\dotsc,n; j\in -N+\mathbb{N}]$ is the ring of global sections of $\W^{n}$ and $(I_{\lambda}^{(N)})$ denotes the ideal generated by the ideal defined in \eqref{eqDefIl}.
The set of such morphisms, viewed as morphisms to $\Lp\W(k)[1/p]^{n}$, is precisely the subset $\diag(p^{\lambda_{1}},\dotsc,p^{\lambda_{n}})\cdot\W(k)^{n} = \Lattice{\lambda}\subset\W(k)[1/p]^{n}$, as claimed.
Finally we discuss the fiber $\mathcal{F} = \mathcal{U}_{\lambda}\times_{\Dem{\lambda}}V$ over a general $k$-valued point $V\in\Dem{\lambda}(k)$, with $N$ as above.
 Again, since $k$ is perfect, each $k$-valued point of $\mathcal{F}$ is given by a morphism
\begin{equation*}
\Spec k \to F := \mathbf{\Up^{(N)}}\times_{\Dem{\lambda}}\lbrace V\rbrace \subset \W^{n},
\end{equation*}
and the set of such morphisms induces a lattice in $\W(k)[1/p]^{n}$. The only thing we have to check is that this lattice is special. However, the product of the elementary divisors of the lattice $F(k)\subset\W(k)[1/p]^{n}$ is determined by the codimension of $F\subset\W^{n}$, which is, by flatness of $\mathbf{U}_{\lambda}^{(N)}\to\Dem{\lambda}$, the same as the codimension of $F_{\lambda}\subset\W^{n}$. Hence, the set of morphisms $F(k)$ induces a special lattice, i.e.~a $k$-valued point of the Grassmannian $\mathcal{F}(k)\in\pGrass(k)$.
\end{proof}

\subsection{Construction of a morphism $\Dem{\lambda}\to\pGrass$}\label{subsectMorphism}

In this subsection we extend the result of Proposition \ref{propMorKPoints} and obtain a \emph{morphism} of fpqc-sheaves to the affine Grassmannian for $\Sl_{n}$,
$$
\Dem{\lambda} \to \pGrass,
$$
which induces the map $\Dem{\lambda}(k)\to\pGrass(k)$ constructed there.

Fix $\lambda\in\domcochar(T)\subset \mathbb{Z}^{n}$.
To simplify notation we let $\Up := \mathbf{U}^{(N)}_{\lambda} \to \Dem{\lambda}$ be as defined in the previous subsection, for some positive integer satisfying the inequalities in \eqref{eqIneqN}.
Our main result in Theorem \ref{thmDescent} will not depend on the particular choice of $N$.

We let $\det: \Lpp\Mat_{n} = \W^{n\times n} \to \W$ be the Greenberg realization of the determinant map for $n\times n$-matrices over $\W(k)$ and consider the composition
\begin{equation*}
\Delta: \Up^{n} \subset \W^{n\times n}\times_{\Spec k}\Dem{\lambda} \to \Lpp\Mat_{n} \xrightarrow{\det} \W.
\end{equation*}
We observe that this morphism of $k$-schemes actually factors through the subscheme
\begin{equation*}
\lbrace 0\rbrace\times\dotsb\times\lbrace 0\rbrace\times\Aff_{\Dem{\lambda}}^{\mathbb{N}} \subset \W,
\end{equation*}
with $nN$ leading $0$'s.
We set $Y_{\lambda} = \lbrace 0\rbrace\times\dotsb\times\lbrace 0\rbrace\times\mathbb{G}_{m}\times\Aff^{\mathbb{N}}_{k}$, again with $nN$ leading $0$'s,
where $\mathbb{G}_{m}\subset\Aff_{k}^{1}$ denotes the multiplicative group, and define a $k$-scheme $X_{\lambda}\subset \Up^{n}$ as the fiber product
\begin{equation}\label{constrX}
\hspace{\fill}
\begin{xy}
\xymatrix{
X_{\lambda} \ar[rr]\ar@^{(->}[d] & & Y_{\lambda} \ar@^{(->}[d] \\
\Up^{n} \ar_{\Delta}[rr] & & \W.
}
\end{xy}
\hspace{\fill}
\end{equation}

\begin{lemma}\label{lemFpqcCover}
 The scheme $X_{\lambda}$ is an open subscheme of $\Up^{n}$ and thus flat over $\Dem{\lambda}$. Moreover $X_{\lambda}$ maps surjectively to the base $\Dem{\lambda}$. In other words, $X_{\lambda}\to\Dem{\lambda}$ is a faithfully flat morphism of $k$-schemes.
\end{lemma}

\begin{proof}
To prove flatness we just note that open immersions are flat, and $\Up^{n}$ is flat over $\Dem{\lambda}$ by construction.
To prove surjectivity consider any geometric point $V\in\OCell{\lambda}$ and let $\kappa$ be its (algebraically closed) residue field.
Then the $\kappa$-valued points of its fiber $\mathcal{F}\subset \mathcal{U}_{\lambda}$ form a special lattice in $\W(\kappa)[1/p]^{n}$ by \ref{propMorKPoints}, and we can find elements $v_{1},\dotsc,v_{n}$ in $\mathcal{F}(\kappa)$ such that $\det(v_{1},\dotsc,v_{n}) = 1$.
Since $\kappa$ is perfect, there exists a morphism of schemes $\beta: \Spec\kappa\to U_{\lambda}^{n}\times_{\Dem{\lambda}}V$ which induces the tuple $(v_{1},\dotsc,v_{n})\in\mathcal{U}_{\lambda}^{n}(\kappa)$.
As we have $\det(v_{1},\dotsc,v_{n}) = 1$, we conclude from the following representations of the determinant morphism between the $k$-ind-schemes $\Lp\Mat_{n}$ and $\Lp\W(k)[1/p]$,
\begin{equation*}
\begin{xy}
\xymatrix{
&
&
&
U_{\lambda}^{n}\times_{\Dem{\lambda}}V \ar@^{(->}[d] &
& \\
W^{n\times n} \ar[r]^{\cdot p}\ar[d]^{\det} &
\cdots \ar[r]^{\cdot p} &
W^{n\times n} \ar[r]^{\cdot p}\ar[d]^{\det} &
W^{n\times n} \ar[r]^{\cdot p}\ar[d]^{\det} & \cdots\\
W\ar[r]^{\cdot p^{n}} &
\cdots \ar[r]^{\cdot p^{n}} &
W\ar[r]^{\cdot p^{n}} &
W\ar[r]^{\cdot p^{n}} & \cdots,
}
\end{xy}
\end{equation*}
that the composition $\det\circ\beta$ factors through $Y_{\lambda}\subset\W$.
This shows that $\beta$ actually induces a $\kappa$-valued point of $X_{\lambda}$ and thus concludes the prove of surjectivity of $X_{\lambda}\to\Dem{\lambda}$. 
\end{proof}

By definition of $X_{\lambda}$ the morphism of $\Dem{\lambda}$-ind-schemes obtained by composition
$$
X_{\lambda} \to \Up^{n} \to \mathcal{U}_{\lambda}^{n} \to \Lp\Mat_{n}
$$
factors through $\Lp\Gl_{n}$.
Composing $X_{\lambda}\to\Lp\Gl_{n}$ with the morphism of $k$-ind-schemes which replaces the first column of each invertible matrix $A$ with its scalar multiple with $\det(A)^{-1}$, we obtain an $\Lpp\Sl_{n}$-equivariant morphism of $k$-ind-schemes $X_{\lambda} \to \Lp\Sl_{n}$, and hence
$$
\varphi_{\lambda}: X_{\lambda}\to\pGrass.
$$

\begin{theorem}\label{thmDescent}
The $X_{\lambda}$-valued point $\varphi_{\lambda}\in\pGrass(X)$ descends to a $\Dem{\lambda}$-valued point $\pi_{\lambda}: \Dem{\lambda}\to\pGrass$. This morphism is equivariant for the (left-)action of $\Lpp \Sl_{n}$ and sends $\Vl\in\Dem{\lambda}(k)$ to the lattice $\Lattice{\lambda}=\diag(p^{\lambda_{1}},\dotsc,p^{\lambda_{n}})\W(k)^{n}\in\pGrass(k)$. Moreover, the restriction of this morphism to $\OCell{\lambda}$ induces a bijection
$$
\OCell{\lambda}(R) \simeq \SCell{\lambda}(R)
$$
for every reduced $k$-algebra $R$.
\end{theorem}

\begin{proof}
Since $\pGrass$ is an fpqc-sheaf by definition and $X_{\lambda}\to\Dem{\lambda}$ is faithfully flat by Lemma \ref{lemFpqcCover}, we have an exact sequence
$$
\pGrass(\Dem{\lambda}) \hookrightarrow \pGrass(X_{\lambda}) \rightrightarrows \pGrass(X_{\lambda}\times_{\Dem{\lambda}}X_{\lambda}).
$$
We have to check that the compositions of $\varphi_{\lambda}$ with the two projections
$$
X_{\lambda}\times_{\Dem{\lambda}}X_{\lambda} \rightrightarrows X_{\lambda}
$$
conincide.
First observe that $X_{\lambda}\to \Dem{\lambda}$ is an affine morphism, and that the descent problem is Zariski-local on $\Dem{\lambda}$. We may thus replace $\Dem{\lambda}$ by an affine open subset $\Spec R\subset \Dem{\lambda}$, and $X_{\lambda}$ by $\Spec S = \Spec R\times_{\Dem{\lambda}}X$, and check whether the images of the induced morphism $\varphi_{\lambda,S} \in \pGrass(S)$ under $\pGrass(S) \rightrightarrows \pGrass(S\otimes_{R}S)$ coincide. In other words, we have to check that, if
$$
\Phi_{1},\Phi_{2}\in\Lp \Sl_{n}(S\otimes_{R}S)
$$
are the two natural compositions $\Spec (S\otimes_{R}S) \rightrightarrows \Spec S \xrightarrow{\Phi} \Lp \Sl_{n}$, then after a possible faithfully flat base change we have $\Phi_{1}^{-1}\cdot\Phi_{2} \in \Lpp\Sl_{n}(S\otimes_{R}S)$. Let $k\subset\kappa$ be an algebraically closed field extension. By construction, a $\kappa$-valued point of $\Spec (S\otimes_{R}S)$ corresponds to a pair of bases of one and the same lattice (given by the corresponding $\kappa$-valued point of $\Spec R\subset \Dem{\lambda}$). Thus the morphism $\Lp\Sl_{n}(S\otimes_{R}S) \to \Lp\Sl_{n}(\kappa)$ sends $\Phi_{1}^{-1}\cdot\Phi_{2}$ to $\Lpp\Sl_{n}(\kappa)$. This means, that we may decompose $\Phi_{1}^{-1}\cdot\Phi_{2} = \Psi + \Omega$, where $\Psi\in\Lpp\Gl_{n}(S\otimes_{R}S)$ (possibly after a faithfully flat base change) and $\Omega\in\Lp\Mat_{n}(S\otimes_{R}S)$ has in its entries only Witt vectors with nilpotent coefficients.
Since $p$-typical Witt vectors with finitely many nilpotent coefficients are killed by multiplication by sufficiently large $p$-powers and $p$ is invertible in $\W(S\otimes_{R}S)[1/p]$ we may conclude that $\Psi+\Omega$ is in fact in the image of $\Lpp\Sl_{n}(S\otimes_{R}S)$ in $\Lp\Sl_{n}(S\otimes_{R}S)$ $p$ $p$.
This concludes the proof of the first part of the theorem.

It is immediate from the definition of $\varphi_{\lambda}$ that
the induced morphism $\pi_{\lambda}: \Dem{\lambda} \to
\pGrass$ sends the lattice scheme $\Vl$ to the lattice
$\Lattice{\lambda}$. In order to see that $\pi_{\lambda}$
induces a bijection $\OCell{\lambda}(R) \simeq
\SCell{\lambda}(R)$ for each reduced $k$-algebra $R$, we
consider the follwing commutative diagram of fpqc-sheaves,
\begin{equation*}
\hspace{\fill}
\begin{xy}
 \xymatrix{
\Lpp\Sl_{n} \ar@^{(->}^{(\id,\Vl)}[rr]\ar@{=}[d] & & \Lpp\Sl_{n}\times \OCell{\lambda} \ar[r]\ar^{(\id,\pi_{\lambda})}[d] & \OCell{\lambda} \ar^{\pi_{\lambda}}[d] \\
\Lpp\Sl_{n} \ar@^{(->}^{(\id,\Lattice{\lambda})}[rr] & & \Lpp\Sl_{n}\times \SCell{\lambda} \ar[r] & \SCell{\lambda},
}
\end{xy}
\hspace{\fill}
\end{equation*}
where the right hand horizontal maps are the morphisms defining
the respective left actions on $\OCell{\lambda}$ and
$\SCell{\lambda}$. Both horizontal
compositions are surjective morphisms of fpqc-sheaves. In order to
check that $\pi_{\lambda}(R): \OCell{\lambda}(R) \to
\SCell{\lambda}(R)$ is bijective, it suffices to show
that the stabilizers, i.e. the respective preimages of
$V_{\lambda, R}$ and
$\Lattice{\lambda}\otimes_{\W(k)}\W(R)$, in $\Lpp\Sl_{n}(R)$, are equal.
Since $R$ is reduced, we may check this fiberwise, i.e.\
assume that $R=\kappa$ is an algebraically closed field.
Thus consider $A \in \Lpp\Sl_{n}(\kappa) = \Sl_{n}(\W(\kappa))$.
By Proposition \ref{propMorKPoints} the lattice $\Lattice{\lambda}\otimes\W(\kappa)$ is identical to the set of $\kappa$-valued points of the fiber $\mathcal{F}_{\lambda}$ in $\mathcal{U}_{\lambda}$  over $\Vl$, and each $\kappa$-valued point of this fiber is defined by a morphism of $\kappa$-valued points of the fiber $F_{\lambda} \in \mathbf{U}_{\lambda}^{(N)}$ over $\Vl$.
Thus $A$ stabilizes $\mathcal{F}_{\lambda}(\kappa)$ if and only if it stabilizes $F_{\lambda}(\kappa)$.
But since $F_{\lambda}$ is a reduced $\kappa$-scheme which is defined by the same equations as $\Vl$ in \eqref{eqDefnIl}, $A$ stabilizes $\Vl$ if and only if it stabilizes the lattice $\Lattice{\lambda}\otimes\W(\kappa)$, which concludes the proof.
\end{proof}

\subsection{Properties of the morphism $\Dem{\lambda}\to\pGrass$}\label{subsectPropertiesOfMorphism}

The restriction of $\pi_{\lambda}$ to $\OCell{\lambda} \to \SCell{\lambda}$ is not an isomorphism of fpqc-sheaves on the full category of $k$-algebras, as Eike Lau pointed out. This follows from

\begin{proposition}[Communicated by Eike Lau]
The functor $\pGrass$ has trivial tangent spaces, and in particular none of the morphisms $\pi_{\lambda}: \OCell{\lambda}\to\SCell{\lambda}$ is an isomorphism except if $\lambda=0$.
\end{proposition}

\begin{proof}
 Let $R=k[\epsilon]/(\epsilon^{2})$. We have to prove that $\pGrass(R) = \pGrass(k)$. First, from the canonical ring homomorphisms $k\to R\to k$ we obtain a factorization of the identity map
 $$
 \pGrass(k) \to \pGrass(R) \to \pGrass(k).
 $$
 Thus we have to check that the right hand map is injective. To this end, consider any two points $L,M\in\pGrass(R)$ which map to the same image $\bar{L}=\bar{M}\in\pGrass(k)$. Now chose a suitable faithfully flat ring extension $R\to S$ such that the points $L_{S}, M_{S}\in\pGrass(S)$, induced by $L$ and $M$, lie in the image of $\Lp\Sl_{n}(S)$. The situation is summarized in the following diagram:
\begin{equation*} 
\hspace{\fill}
\begin{xy}
 \xymatrix{
L, M \ar@{|->}[d]  & \pGrass(R) \ar[r]\ar@^{(->}[d] & \pGrass(k) \ar@^{(->}[d] \\
L_{S},M_{S} & \pGrass(S) \ar[r] & \pGrass(S/\epsilon S) \\
\hat{L},\hat{M} \ar@{|->}[u] & \Lp\Sl_{n}(S) \ar[r]\ar[u] & \Lp\Sl_{n}(S/\epsilon S). \ar[u] 
}
\end{xy}
\hspace{\fill}
\end{equation*}
In fact, the map in the lower row in this diagram is the identity map. Namely, $\Lp\Sl_{n}$ is a $k$-ind-scheme whose connecting homomorphisms annihilate nilpotents, thus in particular any two elements of $\Lp\Sl_{n}(S)$, viewed as matrices over $\W(S)[1/p]$, whose entries differ by multiples of $\epsilon\in S$, are in fact equal. Hence our assumptions on $L$ and $M$ imply that $\hat{L}$ and $\hat{M}$ coincide up to some element $A\in\Lpp\Sl_{n}(S/\epsilon S)$ (possibly after substituting $R\to S$ by a further faithfully flat extension $R\to S\to S'$). Since $\Lp\Sl_{n}(S)=\Lp\Sl_{n}(S/\epsilon S)$ the matrix $A$ trivially lifts to $\hat{A}\in\Lp\Sl_{n}(S)$ such that $\hat{M} = \hat{L}\cdot\hat{A}$. This proves that $L_{S}=M_{S}$, and thus $L=M$, which concludes the prove of injectivity of $\pGrass(R)\to \pGrass(k)$.
\end{proof}

We remind the reader that the Bruhat-order $\leq$ on $\domcochar(T)\subset \mathbb{Z}^{n}$ is defined so that $\mu \leq \lambda$ if and only if the inequality $\mu_{1} +\dotsb +\mu_{i} \leq \lambda_{1}+\dotsb+\lambda_{i}$ holds for every $1\leq i \leq n$. For $\mu\leq\lambda$ and $\mu\neq\lambda$ we write $\mu<\lambda$.

\begin{lemma}\label{lemHF}
Let $\mu, \lambda\in\domcochar(T)$ and let $h_{\mu}$ and $h_{\lambda}$ be the respective Hilbert functions of the associated lattice schemes $V_{\mu}$ and $\Vl$. If $\mu<\lambda$, then $h_{\mu}\neq h_{\lambda}$ and for all $n\in\mathbb{N}$ we have $h_{\mu}(n) \leq h_{\lambda}(n)$ (in this case we write $h_{\mu}<h_{\lambda}$ for short).
\end{lemma}

\begin{proof}
Any cocharacter $\mu < \lambda$ is obtained as a sum of $\lambda$ and cocharacters of the form $(0,\dotsc,0,-1,0,\dotsc,0,1,0,\dotsc,0)$.
Without loss of generality we may assume that $\mu = \lambda + (-1,1,0,\dotsc,0)$.
This means, in the set of generators \eqref{eqDefnIl} for the defining ideal of $\Vl$, we have to replace the generator $x_{1,\lambda_{1}-1}$ by the generator $x_{2,\lambda_{2}}$ in order to obtain a set of generators for the defining ideal of $V_{\mu}$.
As we assumed that $\mu$ is dominant, $\lambda_{2}<\lambda_{1}-1$ and thus $\deg x_{2,\lambda_{2}}^{p^{m}} = \deg x_{1,\lambda_{1}-1}$ with $m=\lambda_{1}-\lambda_{2}-1\geq 1$.
This means that, if we replace in \eqref{eqDefnIl} the generator $x_{1,\lambda_{1}-1}$ by $x_{2,\lambda_{2}}^{p^{m}}$, the result is a generating set of an ideal having \emph{the same} Hilbert function $h_{\lambda}$ as $\Vl$. On the other hand, if we replace $x_{1,\lambda_{1}-1}$ by $x_{2,\lambda_{2}}$ as discussed above, we observe that $x_{2,\lambda_{2}}$ has a strictly lower degree than $x_{2,\lambda_{2}}^{p^{m}}$, and we obtain a generating set of an ideal with a Hilbert function which is $\emph{different}$ from $h_{\lambda}$ and less or equal to $h_{\lambda}$ in each degree. As argued above, this ideal is the ideal of $V_{\mu}$, and the lemma is proved.
\end{proof}

The following two results deal with the set of $k$-valued points of $\Dem{\lambda}$.

\begin{theorem}
At the level of $k$-valued points the morphism $\pi_{\lambda}$ is given by
\begin{align*}
 \pi_{\lambda}(k): \Dem{\lambda}(k) &\to \pGrass(k)\\
		      V &\mapsto (\mathcal{U}_{\lambda}\times_{\Dem{\lambda}}V)(k).
\end{align*}
In other words, at the level of $k$-valued points $\pi_{\lambda}$ coincides with the map constructed in Proposition \ref{propMorKPoints}.
The image of this map is equal to the disjoint union $\coprod_{\lambda'\leq \lambda}\mathcal{C}_{\lambda'}(k)$,
where $\leq$ denotes the Bruhat order on $\domcochar(T)$.
\end{theorem}

\begin{proof}
 To prove the first claim, we just note that a $k$-valued point of the fiber of $X_{\lambda}\to\Dem{\lambda}$ over $V\in\Dem{\lambda}$ is nothing but a basis for the lattice $L$ associated to $V$ by Proposition \ref{propMorKPoints}.
 By construction, $\pi_{\lambda}$ maps $V$ to the right-$\Lpp\Sl_{n}(k)$-coset of this basis, which is the same lattice $L$.
 
 To prove the second claim, we fix $\lambda'>\lambda$ and check that $\Lattice{\lambda'}$ does not lie in the image of $\pi_{\lambda}$.
 By Lemma \ref{lemHF} we know that $h_{\lambda'} > h_{\lambda}$. But since $V_{\lambda'}$ is reduced, it has already the smallest possible Hilbert function among those lattice schemes which possibly map to $\Lattice{\lambda'}$. As $\Dem{\lambda}$ contains only lattice schemes with the same Hilbert function as $V_{\lambda}$, we conclude that $V_{\lambda'}$ does not lie in $\Dem{\lambda}(k)$.
In order to prove that indeed each $\SCell{\lambda'}(k)$ with $\lambda'<\lambda$ is in the image of $\pi_{\lambda}(k)$, we use an argument similar to the one given by  Beauville and Laszlo \cite[Prop.~2.6]{beauville-laszlo}. For integers $e>d$ consider the following equation of matrices over $\W(k((t)))[1/p]$,

\begin{equation}\label{matrixClosure}
\left( \begin{array}{cc} 0& t\\ -t^{-1}& t^{-1}p\end{array} \right)\left( \begin{array}{cc} p^{e}& 0\\ 0& p^{d}\end{array} \right)\left( \begin{array}{cc} t^{-1}& 0\\ t^{-1}p^{e-d-1}& t\end{array} \right) = \left( \begin{array}{cc} p^{e-1}& t^{2}p^{d}\\ 0& p^{d+1}\end{array} \right).
\end{equation}

If we assume $e+d=0$, it follows that the right hand matrix gives rise to a lattice scheme $V\in \Dem{(e,d)}(k((t)))$, which corresponds to a $k((t))$-point of $\SCell{(e,d)}$. Since $\Dem{(e,d)}$ is projective, this $k((t))$-valued point extends to a lattice scheme $\bar{V}$ over $k[[z]]$, whose fiber over $t=0$ maps to $\Lattice{(d+1,e-1)}$.
Surjectivity in the case of a general $n$ follows from the observation that each dominant cocharacter $\lambda' < \lambda$ is obtained as a sum of $\lambda$ and cocharacters of the form $(0,\dotsc,0,-1,0,\dotsc,0,1,0,\dotsc,0)$ (cf. the proof of Lemma \ref{lemHF}). Repeated application of the argument above then shows that $\Lattice{\lambda'}$, and thus also $\SCell{\lambda'}$, is in the image of $\Dem{\lambda}$ for every $\lambda'<\lambda$.
\end{proof}

It would be desirable that the morphisms $\pi_{\lambda}: \Dem{\lambda}\to \pGrass$ were injective, at least at the level of $k$-valued points, in order to reasonably identify the varieties $\Dem{\lambda}$ with ``Schubert varieties'' in the $p$-adic setting. 
Unfortunately, this is not the case, for the reason that the lattice schemes parametrized by $\Dem{\lambda}$ may indeed carry infinitesimal structure.
For example, let $n=2$, $\lambda=(1,-1)\in\mathbb{Z}^{2}$ and $N=2$. Then $\Vl^{(2)}\subset \Aff_{k}^{2\times 2}=\Spec k[x_{-1},x_{0},y_{-1},y_{0}]$ is defined by the ideal $\langle x_{-1},x_{0} \rangle$. Further, $\Dem{\lambda}-\OCell{\lambda}$ contains a whole $\mathbb{A}^{1}_{k}$, whose $k$-points $P_{a},a\in k,$ are given by ideals of the form $\langle y_{-1}+ax_{-1},x_{-1}^{p}\rangle$. This affine line, parametrizing \emph{non-reduced} lattice schemes, maps to the standard lattice $\Lattice{\lambda} = \W(k)^{2} \in \pGrass(k)$. The general situation is summarized by the following Theorem, which is analogous to \citet[Cor.~6.8]{kreidl-2010}.

\begin{theorem}
 A lattice scheme $L\in\Dem{\lambda}(k)$ lies in the open orbit $\OCell{\lambda}(k)$ if and only if it is reduced.
\end{theorem}

\begin{proof}
 Obviously, any $L\in\OCell{\lambda}(k)$ is reduced, since $\Vl$ is reduced. 
 On the other hand, if $\lambda'<\lambda\in\domcochar(T)$, then by Lemma \ref{lemHF} we have $h_{\lambda'}<h_{\lambda}$ for the Hilbert functions of $\Vl$ and $V_{\lambda'}$, respectively.
 Thus $V_{\lambda'}$ does not lie in the closure of $\OCell{\lambda}$. But since any element in $\pi_{\lambda}^{-1}(\Lattice{\lambda'})$ is a lattice scheme with Hilbert function equal to that of $\Vl$ and with a set of $k$-valued points equal to that of $V_{\lambda'}$, we conclude that all the lattice schemes in $\pi_{\lambda}^{-1}(\Lattice{\lambda'})$ carry nontrivial infinitesimal structure.
\end{proof}

Unfortunately, these infinitesimal structures cannot be avoided within our current framework.
As soon as we try to represent lattices by points in a Hilbert scheme, we are forced to use a non-standard grading as described in Subsection \ref{subsectionLatticesWittVector}, in order that Witt vector arithmetics is represented by \emph{graded} morphisms of affine schemes so that $\Lpp\Sl_{n}$ can act on the Hilbert scheme.

A similar situation can be constructed in the function field case, where the analogue of $\Dem{\lambda}$ turns out to be, up to Frobenius twists, a Demazure resolution of the respective Schubert variety in the affine Grassmannian \cite[see][]{kreidl-2010}. This suggests to think of $\Dem{\lambda}$ also in the present Witt vector setting as some sort of Demazure resolution of a Schubert variety in $\pGrass$, but we do not know how to formulate this in a precise way.

\section{Appendix: fpqc-sheaves}\label{sectAppendix}

In this appendix we collect some general results on fpqc-sheaves which are used throughout the preceding sections. In particular, we discuss in detail the existence of sheafifications over the fpqc-site in situations relevant for the present paper.

Let $\mathcal{C}$ be the category of schemes. By a presheaf on $\mathcal{C}$ we mean simply a functor on the category of schemes to the category $\catSet$ of sets.

\begin{lemma}\label{lemRestrSheafif}
Let $\mathcal{D}\subset \mathcal{C}$ be an inclusion of sites, such that fiber products in $\mathcal{D}$ are mapped to fiber products in $\mathcal{C}$.
Assume that for every covering $\mathcal{U} = \lbrace U_{i}\to X \rbrace$ in $\mathcal{C}$ of an object $X\in\mathcal{D}$ there exists a refinement $\mathcal{V}=\lbrace V_{i}\to X\rbrace$ of $\mathcal{U}$ with $V_{i}\in\mathcal{D}$ such that $\mathcal{V}$ is also a covering of $X$ in $\mathcal{D}$.
Then restriction of presheaves from $\mathcal{C}$ to $\mathcal{D}$ commutes with sheafification.
In other words, if a presheaf $F:\mathcal{C}\to \textbf{(sets)}$ has a sheafifcation $F^{a}$, then $F^{a}\rvert_{\mathcal{D}}$ is a (the) sheafification of $F\rvert_{\mathcal{D}}$.
\end{lemma}

\begin{proof}
We check that the canonical morphism $F\rvert_{\mathcal{D}} \to (F^{a})\rvert_{\mathcal{D}}$ is a sheafification on $\mathcal{D}$.
Let $X\in \mathcal{D}$ and let $\xi,\eta \in F(X)$ be such that their images in $F^{a}(X)$ coincide.
By definition of sheafification there exists a covering (in $\mathcal{C}$) of $X$ on which $\xi$ and $\eta$ coincide.
But by assumption this covering can be refined in order to obtain a covering of $X$ in $\mathcal{D}$. Of course, $\xi$ and $\eta$ still coincide on this refinement.
On the other hand, every element $\xi \in F^{a}(X)$ can be represented locally (on a covering in $\mathcal{C}$) by sections of $F$.
Refining this covering, we see that $\xi$ can be represented on a covering in $\mathcal{D}$ by sections of $F$.
Thus $F\rvert_{\mathcal{D}} \to (F^{a})\rvert_{\mathcal{D}}$ is indeed a sheafification and induces an isomorphism $(F\rvert_{\mathcal{D}})^{a} \to (F^{a})\rvert_{\mathcal{D}}$.
\end{proof}

\begin{theorem}[\citealt{vistoli-2008}]\label{thmSheafificationExists}
Let $F$ be a presheaf on $\mathcal{C}$. Assume that $F$ is a sheaf for the Zariski topology. Then $F$ is an fpqc-sheaf on $\mathcal{C}$ if and only if for every faithfully flat homomorphism of affine schemes $Y\to X$ the sequence
\begin{equation}\label{eqEqualizer}
F(X) \to F(Y) \rightrightarrows F(Y\times_{X}Y)
\end{equation}
is an equalizer.
\end{theorem}

\begin{proposition}\label{propfpqcSheafification}
Let $F$ be a presheaf on $\mathcal{C}$. Assume that $F$ satisfies the following two conditions:
\begin{enumerate}[(1)]
\item for every faithfully flat morphism of \emph{affine} schemes $Y\to X$ the sequence
$$
F(X) \to F(Y) \rightrightarrows F(Y\times_{X}Y)
$$
is an equalizer, and
\item for every finite collection of affine schemes $Y_{1},\dotsc,Y_{n}$ we have
$$
F(Y_{1}\coprod \dotsb \coprod Y_{n}) = F(Y_{1})\times\dotsb \times F(Y_{n}).
$$
\end{enumerate}
Then the Zariski-sheafification $F^{a}$ of $F$ is an fpqc-sheaf. In particular, $F^{a}$ is an fpqc-sheafification of $F$. Moreover, the natural transformation $F\to F^{a}$ restricts to an isomorphism on the category of affine schemes.
\end{proposition}

\begin{proof}
In view of Theorem \ref{thmSheafificationExists} we only have to prove that the condition in (1) of the present proposition remains valid after Zariski-sheafification. Thus it will suffice to prove the last assertion, namely that the natural map $F(X)\to F^{a}(X)$ is indeed an isomorphism for every affine $X$. To this end, for an arbitrary scheme $X$ and any Zariski-covering $\mathcal{U}$ of $X$ let $K(\mathcal{U})$ be the equalizer of $F(\mathcal{U})\rightrightarrows F(\mathcal{U}\times_{X}\mathcal{U})$. If we set $F'(X) = \dlim_{\mathcal{U}}K(\mathcal{U})$, where the colimit is taken over all Zariski-coverings of $X$, then $F'$ will be a separated presheaf. Applying this procedure twice, i.e. forming $F''$, will yield a \emph{sheaf}, and indeed $F''$ is equal to the Zariski-sheafification $F^{a}$ of $F$. Now observe the following: if $X$ is affine, there is a cofinal subsystem of all Zariski coverings of $X$ given by those coverings which consist of only \emph{finitely many affines}. Thus, using assumption 
(2),
$$
F'(X) = \dlim_{Y\to X} \ker(F(Y)\rightrightarrows F(Y\times_{X}Y)),
$$
where now the limit is taken over a certain family of faithfully flat morphisms $Y\to X$ of affine schemes. But by assumption (1) for every such $Y\to X$ we have $F(X) = \ker(F(Y)\rightrightarrows F(Y\times_{X}Y))$, whence $F'(X)=F(X)$. This implies $F^{a}(X)=F(X)$, as desired.
\end{proof}

\begin{corollary}
Let $F$ be as in Proposition \ref{propfpqcSheafification}. Then the restriction of $F$ to the site of affine schemes (with arbitrary covering families consisting of affine schemes) is a sheaf for the fpqc-topology.
\end{corollary}

The preceding discussion shows that the category of fpqc-sheaves on the category of $k$-schemes is equivalent to the category of functors on affine $k$-schemes which satisfy the conditions (1) and (2) of Proposition \ref{propfpqcSheafification}. Mutually inverse equivalences are given by restriction and respectively by passing to the associated Zariski-sheaf. In \cite{beauville-laszlo}, the authors indeed define a $k$-space to be a functor on the category of affine $k$-schemes which satisfies condition (1). On the other hand, they do not require condition (2), which, however, does not seem to be automatic.

The following proposition shows that indeed every directed system of $k$-schemes gives rise to a $k$-ind-scheme (i.e. the colimit in the category of $k$-spaces exists).

\begin{proposition}\label{propSheafifOfDirectSystem}
A functor which is defined as an inductive limit of schemes admits an fpqc-sheafification. More precisely, its Zariski-sheafification is already an fpqc-sheaf(ification). Further, the restriction of this sheafification to the category of affine schemes coincides with the original presheaf defined by the inductive system of schemes.
\end{proposition}

\begin{proof}
We have to check that such a functor satisfies the assumptions (1) and (2) of Proposition \ref{propfpqcSheafification}.

To this end, let $(X_{i})$ be a direct system of schemes and let $\dlim X_{i}$ be its colimit in the category of presheaves. Let $T_{1},\dotsc,T_{n}$ be affine schemes. Then we have
\begin{multline*}
(\dlim X_{i})(T_{1}\coprod\dotsb\coprod T_{n}) = \dlim (X_{i}(T_{1}\coprod\dotsb\coprod T_{n})) = \\ = \dlim (X_{i}(T_{1})\times\dotsb\times X_{i}(T_{n})) = (\dlim X_{i})(T_{1})\times\dotsb\times (\dlim X_{i})(T_{n}),
\end{multline*}
which is condition (2). It remains to check exactness of the sequence
$$
(\dlim X_{i})(R) \to (\dlim X_{i})(S) \rightrightarrows (\dlim X_{i})(S\otimes_{R}S),
$$
where $R\to S$ is a faithfully flat homomorphism of rings. Thus let $P\in (\dlim X_{i})(S)$ be such that both images of $P$ in $(\dlim X_{i})(S\otimes_{R}S)$ coincide. Assume that $P$ is represented by an element $P'\in X_{i}(S)$. By definition of the inductive limit, there exists some $i\leq j\in I$ such that that the induced objects in $X_{j}(S\otimes_{R}S)$ coincide. Now we can use the exactness of the sequence
$$
X_{j}(R) \to X_{j}(S) \rightrightarrows X_{j}(S\otimes_{R}S)
$$
to obtain an $R$-valued point of $X_{j}$, and hence an $R$-valued point of $\dlim X_{i}$ which induces $P$. Injectivity of the map $(\dlim X_{i})(R) \to (\dlim X_{i})(S)$ is proved likewise, which shows that condition (1) holds as well.
\end{proof}

Proposition \ref{propSheafifOfDirectSystem} says that if we restrict the functor direct-limit $\dlim X_{i}$ to the category of affine schemes (or more generally: quasi-compact schemes), then it is already a sheaf for the fpqc-topology. This is Beauville and Laszlo's point of view.\\

Once again let $k$ denote a field of positive characteristic $p$. Contrary to what Vistoli claims \cite[][Thm.~2.64]{vistoli-2008}, arbitrary functors on the category of $k$-schemes do not in general admit an fpqc-sheafification.
An example of such a functor is described by \citet{waterhouse}.
As Waterhouse explains, the general problem with constructing an fpqc-sheafification of an arbitrary functor is that one is forced to consider direct limits over ``all'' fpqc-coverings of a given scheme.
However, the entirety of ``all'' fpqc-coverings will not be a set, but a proper class.
One way out of this problem would be to restrict to a fixed universe, which has the drawback that sheafifications depend on the particular choice of the universe.
On the other hand, Waterhouse proves that for \emph{basically bounded} functors it suffices to look at direct limits over certain \emph{sets} of fpqc-coverings, which resolves the above described set-theoretical problems.


Let $m$ be a cardinal number not less than the cardinality of $k$, fix a set $S$ of cardinality $m$, and let $\catAlgm{k}$ be the category of $k$-algebras whose underlying set is contained in $S$. Let $\catAlg{k}$ denote the category of $k$-algebras, and let $j:\catAlgm{k} \hookrightarrow \catAlg{k}$ be the inclusion. For any set-valued functor on the category of $k$-algebras, let $j^{*}$ denote the restriction to $\catAlgm{k}$. Right-adjoint to $j^{*}$ is the Kan extension $j_{*}$ along $\catAlgm{k} \hookrightarrow \catAlg{k}$.

\begin{definition}
A functor $F$ on the category of $k$-algebras is \emph{$m$-based} if it has the form $j_{*}G$ for some functor $G$ on $\catAlgm{k}$. A functor is \emph{basically bounded} if there exists a cardinal $m$ such that it is $m$-based.
\end{definition}

\begin{theorem}[\citealt{waterhouse}, Cor.~5.2]\label{thmWaterhouse}
If a functor $F$ on the category of $k$-algebras is $m$-based, then it has an fpqc-sheafification. More precisely, if $j^{*}F\to G$ is a sheafification for the fpqc-topology on $\catAlgm{k}$, then $F=j_{*}j^{*}F \to j_{*}G$ is an fpqc-sheafification on $\catAlg{k}$.
\end{theorem}

We use the following two observations made by Waterhouse: (a) A functor which is represented by an affine scheme whose underlying ring has cardinality $\leq m$ is $m$-based, and (b) the Kan extension $j_{*}$ preserves colimits, and in particular, the colimit over a system of basically bounded functors is again basically bounded.

\begin{corollary}\label{corGroupQuotientBasicallyBounded}
Let $H$ and $G$ be functors on the category of $k$-algebras with values in groups, which are represented by $k$-ind schemes. Let $H\to G$ be a functorial group homomorphism. Then the presheaf-quotient $G/H$ is basically bounded, and hence has a well-defined fpqc-sheafification.
\end{corollary}

\begin{proof}
By (a) the ind-schemes $H$ and $G$ are colimits of basically bounded functors. Thus they are themselves basically bounded by (b) above. Further, since $G/H$ is the colimit of a direct system of the form
$$
G \times H \rightrightarrows G,
$$
this presheaf-quotient is basically bounded, again by (b). By Waterhouse's theorem, it thus has a well-defined fpqc-sheafification on the category of $k$-algebras. 
\end{proof}

\section{Acknowledgments} I am very much indebted to Ulrich
G\"ortz at the Institut f\"ur Experimentelle Mathematik at
the University of Duisburg-Essen for his interest in my work
and his patient advice during the preparation of the present
paper. Further, my thanks are due to Moritz Kerz, who
raised the question answered by Theorem \ref{thmIntro2}, and
to Philipp Hartwig, who pointed out to me the article by
Waterhouse on fpqc-sheafifications.
Also I would like to thank Eike Lau, who discovered an error in a previous version of Theorem \ref{thmIntro1}, and the referee for many additional suggestions and comments.

This work was partially supported by the SFB/TR 45 ``Periods,
Moduli Spaces and Arithmetic of Algebraic Varieties'' of the
DFG (German Research Foundation) and is part of the author's PhD-thesis.

\end{document}